\documentclass[final]{article}

\usepackage{graphicx}
\usepackage{amsmath}

\newtheorem{thm}{Theorem}[section]

\newcommand{\authorone}[2][]{\hspace*{9pt}{\small\textrm{\uppercase{#2}},$^{*}$ \textit{#1}}\par}
\newcommand{\authortwo}[2][]{\hspace*{9pt}{\small\textrm{\uppercase{#2}},$^{**}$ \textit{#1}}\par}
\newcommand{\authorthree}[2][]{\hspace*{9pt}{\small\textrm{\uppercase{#2}},$^{***}$ \textit{#1}}\par}

\newcommand{\addressone}[1]{\footnote{\hspace*{-14pt}$^{*}\,$Postal address: #1}\par}
\newcommand{\addresstwo}[1]{\footnote{\hspace*{-14pt}$^{**}\,$Postal address: #1}\par}
\newcommand{\addressthree}[1]{\footnote{\hspace*{-14pt}$^{***}\,$Postal address: #1}\par}

\renewcommand{\author}{\authorone}

\newcommand{\keywords}[1]
           {\begin{center}
            \begin{minipage}{315.83pt}
            \small
            \noindent \emph{Keywords:}~{\textrm{#1}}
            \end{minipage}
            \end{center}
            \normalsize
           }

\newcommand{\ams}[2]
           {\begin{center}
            \begin{minipage}{315.83pt}
            \small
            \noindent 2000 Mathematics Subject Classification:~Primary {\uppercase{#1}}\\
            \phantom{2000 Mathematics Subject Classification:~}Secondary {\uppercase{#2}}
            \end{minipage}
            \end{center}
            \par\normalsize
           }

\begin{document}

\title{Fluid limit theorems for stochastic hybrid systems with application to neuron models}

\date{4th May, 2009}
\maketitle



\authorone[Institut Jacques Monod UMR 7592 CNRS, Univ. Paris VII, Univ. Paris VI]{K. Pakdaman} 
\addressone{Institut Jacques Monod UMR7592\\CNRS, Univ. Paris VII, Univ. Paris VI B\^atiment Buffon
15 rue H\'el\`ene Brion
75205 Paris cedex 13 Paris, France} 
\authortwo[Laboratoire de Probabilit\'es et Mod\`eles Al\'eatoires UMR7599
Univ. Paris VI, Univ VII - CNRS]{M. Thieullen} 
\addresstwo{Laboratoire de Probabilit\'es et Mod\`eles Al\'eatoires UMR7599
Univ. Paris VI, Univ VII - CNRS \\Bo\^ite 188\\ University Paris VI\\ 175, rue du Chevaleret 75013 Paris, France} 
\authorthree[CREA, Ecole Polytechnique, Paris, France; Institut Jacques Monod UMR 7592 CNRS, Univ. Paris VII, Univ. Paris VI; Laboratoire de Probabilit\'es et Mod\`eles Al\'eatoires UMR7599
Univ. Paris VI, Univ VII - CNRS]{G. Wainrib} 
\addressthree{Laboratoire de Probabilit\'es et Mod\`eles Al\'eatoires UMR7599
Univ. Paris VI, Univ VII - CNRS \\Bo\^ite 188\\ University Paris VI\\ 175, rue du Chevaleret 75013 Paris, France} 


\begin{abstract}
This paper establishes limit theorems for a class of stochastic
hybrid systems (continuous deterministic dynamic coupled with jump
Markov processes) in the fluid limit (small jumps at high
frequency), thus extending known results for jump Markov processes.
We prove a functional law of large numbers with exponential
convergence speed, derive a diffusion approximation and establish a
functional central limit theorem. We apply these results to neuron
models with stochastic ion channels, as the number of channels goes
to infinity, estimating the convergence to the deterministic model. In terms of neural coding, we apply our central limit
theorems to estimate numerically impact of channel noise both on frequency and spike timing coding.
\end{abstract}

\ams{60F05;60F17;60J75}{92C20;92C45}

\section{Introduction}     

\noindent In this paper we consider stochastic hybrid systems where a
continuous deterministic dynamic is coupled with a jump Markov
process. Such systems were introduced in \cite{Davis1984}  as
piecewise deterministic Markov processes. They have been
subsequently generalized to cover a wide range of applications: communication networks,
biochemistry and more recently DNA replication modeling \cite{blom2006,hespanha2005,Kouretas2006,LygPNAS}. We are interested in the fluid limit for these systems considering the case
of small jumps of size $1/N$ at high frequency $N$, with a view
towards application to neural modeling.


\noindent The general class of model we consider is described in section 2.1, and for the sake of clarity, we describe here a simple example which retains the main features. Consider a population of $N$ independent individuals, each of them being described by a jump Markov process $u_k(t)$ for $k=1,...,N$ with states $0$ and $1$, and with identical transition rates $\alpha>0,\beta>0$ as follows:
\begin{figure*}[!h]
\center
\includegraphics[scale=5]{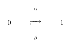}
\end{figure*}
As an empirical measure, we define the proportion of individuals in state $1$ at time $t$ by:
$$e_N(t)=\frac{1}{N}\displaystyle{\sum_{k=1}^N} u_k(t)$$
The model becomes hybrid when we assume a global coupling through a variable $V_N\in \mathbf{R}$, in the sense that the rates $\alpha(V_N)$ and $\beta(V_N)$ are functions of $V_N$. This variable $V_N$ is itself solution of a differential equation, between the jumps of $e_N(t)$:
$$\frac{dV_N}{dt} = f(V_N,u_N)$$
where $f:\mathbf{R}^2\to \mathbf{R}$. In the general case, this model is extended with more general non-autonomous jump Markov processes, the global variable can be vector valued and the transition rates can be functions of the empirical measure (section 2.1). 

\noindent We prove convergence in probability on finite time intervals,with techniques inspired by
\cite{Austin2007},
of the solution $X_N$ of the stochastic hybrid system to a
deterministic limit $x=(v,g)$. For the example above, $x$ is solution of:
\begin{eqnarray*}
\frac{dv}{dt}&=&f(v,g)\\
\frac{dg}{dt}&=&(1-g)\alpha(v)-g\beta(v)
\end{eqnarray*}
We derive a
diffusion approximation and prove a functional central limit theorem 
that helps characterizing the fluctuations of both the discrete and
continuous variables around the deterministic solution. We
obtain that these fluctuations are a gaussian process which
corresponds to the asymptotic law of the linearized diffusion
approximation. We further obtain an exponential speed of
convergence which relates the tail distribution of the error
$E_N(T)=\displaystyle{\sup_{[0,T]}} |X_N-x|^2$ to the size parameter
$N$ and the time window $T$ : for $\Delta>0$ and $N$ large,
\begin{equation}
P(E_N(T)>\Delta)\leq e^{-\Delta N H(T)}
\end{equation}
Thus the convergence
result can be extended to large time intervals $[0,T(N)]$, provided
that $T=T(N)$ is such that $NH(T(N))\to \infty$. Inequality (1.1) is a new result which provides an estimate to the required number $N$ of individuals to reach a
given level of precision. This number increases with the time scale on which one
wants this precision to be achieved. For system subject to
finite-size stochasticity, sometimes called demographic
stochasticity it provides a relation
between the reliability time-scale to the population size $N$. There are other ways of
obtaining a law of large numbers, for example using the convergence
of the master equation or of the generators \cite{ethierkurtz}.
We want to highlight here that our proof is based on exponential inequalities for martingales. Other ways of obtaining a law of large numbers would not be likely to provide an estimate such as $(1.1)$.

\noindent Our mathematical reference on the fluid limit is the seminal paper \cite{Kurtz1971} which contains a law of
large numbers and a central limit theorem  for sequences of jump
Markov processes. Recently, a spatially extended version of these
models has been considered in \cite{Austin2007}, for a standard
neuron model. The author shows convergence in probability up to
finite time windows to a deterministic fluid limit expressed in
terms of a PDE coupled with ODEs. In the present paper, we consider
a class of non-spatial models which however includes multi
compartmental models, by increasing the dimension. We extend the
results of \cite{Kurtz1971} to stochastic hybrid models at the fluid
limit.


\noindent Neurons are subject to various
sources of fluctuations, intrinsic (membrane noise) and extrinsic
(synaptic noise). Clarifying the impact of noise and variability in
the nervous system is an active field of research \cite{segundo94},
\cite{Faisal08}. The intrinsic fluctuations in single neurons are
mainly caused by ion channels, also called channel noise, whose
impacts and putative functions are intensively investigated
\cite{white2000cnn,shuai2003,rowat2007iis}, mainly by numerical
simulations. Our motivation is to study the intrinsic fluctuations
in neuron models and we think that stochastic hybrid systems are a natural tool for this purpose. The channels open and
close, through voltage induced electromagnetic conformational
change, thus enabling ion transfer and action potential generation.
Because of thermal noise, one of the main features of those channels
is their stochastic behavior. 

\noindent In terms of modeling, our starting point is the
stochastic interpretation of the Hodgkin-Huxley formalism
\cite{HH52}. In this setting, ion
channels are usually modeled with independent Markov jump processes,
whose transition rates can be estimated experimentally
\cite{Vandenberg1991}. These stochastic discrete models are coupled with a
continuous dynamic for the membrane potential, leading to a
piecewise-deterministic Markov process. Thus, the
individuals are the ion channels and the global variable $V_N$ the
voltage potential (cf. section 3.).
\noindent 
Deterministic hybrid kinetic equations
appear to be a common formalism suitable for each stage of nervous
system modeling as shown in \cite{destexhe1994}. This latter study provides
us with a framework to introduce stochastic hybrid processes to
model action potential generation and synaptic transmission, as stochastic version of deterministic kinetic models
coupled with differential equations through the transition rates.

\noindent On the side of neuron modeling applications, the limit behavior of a similar but less general model is
considered in \cite{Fox1994}, using an asymptotic development of the
master equation as $N\to \infty$, which formally leads to a
deterministic limit and a Fokker-Planck equation (Langevin
approximation), providing the computation of the diffusion
coefficients. The Langevin approximation is also studied in
\cite{Tuckwell87}, but in a simplified case where the transition
rates  are constants (independent of $V_N$), which is actually the
case studied in \cite{Kurtz1971}. Our mathematical results extend these previous studies to
a wider class of models (if we put aside the spatial aspects in
\cite{Austin2007}), providing a rigorous approach for the Langevin
approximation, and establishing a central limit theorem which describe the effect of channel noise on the
membrane potential \cite{Steinmetz2000}.
The convergence speed provides a quantitative insight into
the following question : if a neuron needs to be reliable during a
given time-scale, what would be a sufficient number of ion channels?
We thus provide a mathematical foundation for the study of stochastic neuron models, and we apply our results to standard
models, quantifying the effect of noise on neural coding. In
particular, both frequency coding (sec. 3.5.1) and spike timing coding (sec. 3.5.2) are numerically studied with Morris-Lecar neuron model with a large number
of stochastic ion channels.

\noindent Generically, stochastic hybrid models in the fluid limit
would arise in multiscale systems with a large population of stochastic agents
coupled, both top-down and bottom-up, through a global variable, leading to an emergent cooperative behavior. Starting from a microscopic
description (ion channels), the central limit theorem as stated in this paper leads to a description of the
fluctuations of the global variable (membrane potential). So, in the perspective of
applications, it would be interesting to investigate how our
framework and results could be developed in other fields than neural modeling: for instance in chemical kinetics, in
population dynamics, in tumor modeling, in economics or in opinion
dynamics theory. In a more mathematical perspective, it would be interesting to consider a wider class of models, for instance by including spatial aspects as in \cite{Austin2007} or by weakening the independence assumption. Other questions could be investigated, for instance concerning escape problems, first passage times and large deviations, whenever $N$ is large or not. 

\noindent Our paper is organized as follows. In section 2. we define
our model and formulate the main results. In section 3., we apply our results to neuron models. In
section 4. we give the proof of the law of large numbers and its convergence speed
(\textbf{Theorem 2.1}) and in section 5. we give the proof of the
Langevin approximation (\textbf{Theorem 2.2}) and central
limit theorems (\textbf{Theorem 2.3-2.4-2.5}).

\section{Model and main results}
This section contains the definition of our general model and states the main theorems.

\subsection{Model}

\paragraph{Stochastic hybrid model $\mathbf{(S_N)}$}
Let $p,q,N \in \mathbf{N}^*$, and $r_j\in \mathbf{N}^*$ for all
$1\leq j\leq q$. Let $d=\sum_{j=1}^q r_j$. We define the stochastic
hybrid model $(S_N)$, whose solution
$$X_N(t)=(V_N(t),\mathbf{e}_N(t))\in\mathbf{R}^p\times
\mathbf{R}^d,\  t\geq 0$$ satisfies:
$$\frac{dV_N}{dt}=f(X_N)$$
and $\mathbf{e}_N=(e_N^{(1)},...,e_N^{(q)})$ with $e_N^{(j)}\in
\mathbf{R}^{r_j}$, where the processes $e_N^{(j)}(t)$ are $q$
independent jump Markov processes. Note that the differential equation for $V_N$ is holding only
between the jump times of the process $\mathbf{e}_N$, with updated
initial conditions. For $1\leq
j\leq q$, processes $e_N^{(j)}(t)$ are characterized by, 
\begin{itemize}
\item their state space : $E_N^{(j)}=\left\{(x_1,...,x_{r_j})\in\{0,\frac{1}{N},...,1\}^{r_j} \ | \ \sum_{k=1}^{r_j}x_k=1\right\}$
\item their intensity $\lambda^{(j)}_N$:
for $X=(V,\mathbf{e})\in\mathbf{R}^p\times \mathbf{R}^d$,
$\lambda_N^{(j)}(X)=N {\tilde{\lambda}^{(j)}(X)}$ with
$$\tilde{\lambda}^{(j)}(X)=\displaystyle{\sum_{k=1}^{r_j}\mathbf{e}^{(j)}_k \sum_{l=1,\ l\neq k}^{r_j}\alpha_{k,l}^{(j)}(X)}$$
\item their jump law $\mu_N^{(j)}$: we define $u^{(j)}_a=(0,...,0,1,0,...,0)\in \mathbf{R}^{r_j}$ and $u^{(j)}_{a,b}=u^{(j)}_a-u^{(j)}_b$ for $1\leq a,b\leq r_j$. The transition of an individual agent in the population $j$ from
one state $a$ to another state $b$ corresponds to a jump of
$z=\frac{1}{N}u^{(j)}_{b,a}$ for the process $e_N^{(j)}$. Thus we
define:
$$X+\frac{1}{N}\Delta X^j_{a,b}=(V,e^{(1)},...,\tilde{e}^{(j)},...,e^{(q)})$$
$$\tilde{e}^{(j)}=e^{(j)}+u^{(j)}_{a,b}$$
So that the jump law for a jump of $z$ is given by:
$$\mu_N^{(j)}(X,z)= \frac{\mathbf{e}_a^{(j)}\alpha_{a,b}(X)}{\tilde{\lambda}^{(j)}(X)} \mbox{ if } z=\frac{1}{N}u^{(j)}_{b,a},$$ for all $1\leq a,b \leq r_j$ such that $\mathbf{e}_a^{(j)}\neq 0$ and $\mathbf{e}_b^{(j)}\neq 1$, and $$\mu_N^{(j)}(X,z)=0 \mbox{
otherwise.}$$
\end{itemize} 
 For a more formal definition we refer to
\cite{Davis1984}. 

\noindent For $1\leq k \leq r_j$, the $k$-th component $\{e_N^{(j)}\}_k$ of vector $e_N^{(j)}$ can be interpreted as the
proportion of agents of type $j$ which are in the state $k$ in a
population of size $N$.

\noindent We show below in \textbf{Theorem 2.1} that this stochastic hybrid model has a limit as $N \to \infty$ which is the following deterministic model.
\paragraph{Deterministic model $\mathbf{(D)}$}
We define the deterministic model $(D)$, whose solution
$X=(v,\mathbf{g})\in\mathbf{R}^p\times \mathbf{R^d}$ with
$\mathbf{g}=(\mathbf{g}^{(1)},..,\mathbf{g}^{(q)})$ satisfies:

$$\left\{
    \begin{array}{ll}

         \dot{v}=f(v,\mathbf{g})\\
     \dot{\mathbf{g}}^{(j)}_k=\displaystyle{\sum_{1\leq i \leq r_j,i\neq k} \alpha^{(j)}_{i,k}(X)\mathbf{g}^{(j)}_i-\alpha^{(j)}_{k,i}(X)\mathbf{g}^{(j)}_k}     \end{array}
\right.  \  \  \  \  (D)$$ for all $1 \leq j \leq q,\ \forall 1\leq
k \leq r_j$. The first equation is the same as in the stochastic
model (deterministic part) and the second equation corresponds to
the usual rate equation, with a gain term and a loss term.
\\

\noindent The following example illustrates the general model in a simpler
relevant setting motivated by applications. This setting will be used 
in the proofs in section 4 and 5 in order to make the arguments clearer.

\paragraph{\textbf{Example}}
We consider the case where $p=q=1$ and $r_1=2$. We can construct a
stochastic hybrid process as follows: first let us introduce a collection of $N$ independent jump Markov
processes $u^{(k)}$ for $1\leq k\leq N$ with $u_t^{(k)}:0\to 1$ with
rate $ \alpha(V_N)$ and $ 1\to 0$ with rate $\beta(V_N)$:
\begin{figure*}[!h]
\center
\includegraphics[scale=5]{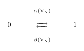}
\end{figure*}

\noindent where $V_N$ is defined below. We then consider $e_N(t)=(\{e_N\}_0(t),\{e_N\}_1(t))$ the proportions of
processes in the states $0$ and $1$. In this case, the stochastic
hybrid model $(S_N)$ can be written as:

$$\left\{
    \begin{array}{lll}

\dot{V}_N(t)=f(V_N(t),e_N(t))\\
e_N(t)=\left(\frac{1}{N}\sum_{k=1}^N
\delta_0(u_t^{(k)}),\frac{1}{N}\sum_{k=1}^N
\delta_1(u_t^{(k)})\right)\\
V_N(0)=v_0 ;\  e_N(0)=(u_0,1-u_0)

    \end{array}
\right.$$

\noindent Note that, if we define $u_N(t)=\frac{1}{N}\sum_{k=1}^N
\delta_1(u_t^{(k)})$, then we have $e_N(t)=\left(1-u_N(t),
u_N(t)\right)$, so that the solution is determined by the pair
$X_N(t)=(V_N(t),u_N(t))$.

\noindent Thus, each member of the sequence of jump Markov processes
$\{u_N\}_{N\geq 1}$ is characterized by
\begin{itemize}
\item its state space $E_N=\{0,\frac{1}{N},\frac{2}{N},...,1\}$,
\item its intensity $\lambda_N(V_N(t),u)=N[u\beta(V_N(t)) + (1-u)\alpha(V_N(t))]$. This intensity is time-dependent through $V_N(t)$.
\item its jump law $$\mu_N(V_N(t),u,y)=\mu^+(V_N(t),u)\delta_{y,u+\frac{1}{N}}+\mu^-(V_N(t),u)\delta_{y,u-\frac{1}{N}}$$
where $\mu^+(V,u)=\frac{(1-u)\alpha(V)}{u\beta(V) + (1-u)\alpha(V)}$
and $\mu^-(V,u)=\frac{u\beta(V)}{u\beta(V) + (1-u)\alpha(V)}$. This
jump law is also time-dependent through $V_N(t)$.
\end{itemize}

\noindent The deterministic system $(D)$takes the form:

$$\left\{
    \begin{array}{lll}
        \dot{v}(t)=f(v(t),g(t))\\
        \dot{g}(t)=(1-g(t))\alpha(v(t))-g(t)\beta(v(t))\\
     v(0)=v_0\ ;\ g(0)=u_0
    \end{array}
\right.$$

\noindent In the sequel, we will be interested in the asymptotic behavior of the stochastic hybrid models $(S_N)$ under the limit fluid assumption. Let us now recall what this assumption means. Let $(x_N)$ be a sequence of homogeneous Markov jump
processes with state spaces $E_N \subset R^k$,
intensities $\lambda_N(x)$ and jump law
$\mu_N(x,dy)$. Define the flow as
$F_N(x)=\lambda_N(x)\int_{E_N}(z-x)\mu_N(x,dz)$. The fluid limit
occurs if the flow admits a limit and if the second order moment of
the jump size converges to zero when $N\to \infty$. Our stochastic hybrid model is in the fluid limit since the jumps
are of size $1/N$ and the intensity is proportional to $N$. This
stems from the fact that we are modeling proportions in a population
of independent agents. However, this independence assumption is not
necessary to satisfy the fluid limit.

\subsection{Law of large numbers for stochastic hybrid systems}
We give here the first result concerning the convergence of the
stochastic hybrid model $(S_N)$, which is a functional law of large
numbers on finite time windows.

\begin{thm}[]
Let $\epsilon >0$, $\delta >0$, $T>0$.
Let us assume that the functions $\alpha_{i,j}$ and $f$ are $C^1$, and satisfy the following condition:
\begin{center}
(H1): the solution $v$ of $(D)$ is bounded on $[0,T]$ and for all $N\geq 1$ the solution $V_N(t)$ of $(S_N)$ is uniformly bounded in $N$ on
$[0,T]$.
\end{center}

Let $X^{init}$ a given initial condition for $(D)$ and
$X=(v,\mathbf{g})$ its solution. Then there exists an initial condition
$X_N^{init}$ for $(S_N)$ and $N_0\geq 0$ such
that $\forall N\geq N_0$, the solution
$X_N=(V_N,\mathbf{e}^{(1)}_N,..,\mathbf{e}^{(q)}_N)$ satisfies, for
all $1\leq j\leq q$ and $1\leq k \leq r_j$ :
$$\mathbf{P}\left[\displaystyle{\sup_{0\leq t \leq T} ||V_N(t)-v(t)|| } \geq \delta\right] \leq \epsilon$$
$$\mathbf{P}\left[\displaystyle{\sup_{0\leq t \leq T} |\{\mathbf{e}^{(j)}_N\}_k(t)-\mathbf{g}^{(j)}_k(t)| } \geq \delta\right] \leq \epsilon$$
Moreover, if we define
$$P_N(T,\Delta)=\mathbf{P}\left[\displaystyle{\sup_{0\leq t \leq T} ||V_N(t)-v(t)||^2+ \sum_{j=1}^q ||\mathbf{e}^{(j)}_N(t)-\mathbf{g}^{(j)}(t)||^2 } > \Delta\right]$$
there exist two constants $B(T)>0$ and $C>0$ such that for $\Delta$
sufficiently small:
\begin{equation}
\displaystyle \limsup_{N \to \infty} \frac{1}{N}\log P_N(T,\Delta) \leq - \frac{\Delta e^{-B(T)T}}{CT}
\end{equation}
\noindent Moreover if
\begin{center}
(H2): assumption (H1) holds true on $[0,+\infty[$\\
\end{center}
then the constant $B(T)=BT$ is proportional to $T$.
\end{thm}

\paragraph{Interpretation of the convergence speed}
We have obtained in $(2.1)$ an upper bound for the convergence speed which can
help to answer the following issue. Given a number of channels
$N$, given an error $\Delta$ and a confidence probability $1-p$ (e.g
$p=0.01$), the time window $[0,T]$ for which we can be sure
(up to probability $1-p$) that the distance between the stochastic
and the deterministic solutions (starting at the same point) is less
than $\Delta$ is given by $(2.1)$. In section 3.3, we show numerical simulation
results illustrating the obtained bound for the convergence speed
for the stochastic Hodgkin-Huxley model.

\paragraph{Remark} Assumption $(H2)$ and thus $(H1)$, are satisfied for most neuron models, for instance for the Hodgkin-Huxley (HH) model \cite{cronin}.

\subsection{Langevin approximation}

Our second result is a central limit theorem that provides a way to
build a diffusion or Langevin approximation of the solution of the stochastic hybrid
system $(S_N)$. For $X=(v,\mathbf{e})\in \mathbf{R}^p\times \mathbf{R^d}$, for
$1\leq j \leq q$, $1\leq i,k\leq r_j$, let

 $$b_{j,k}(X)=\displaystyle{\sum_{1\leq i \leq r_j,i\neq k} \alpha^{(j)}_{i,k}(
X)\mathbf{e}^{(j)}_i-\alpha^{(j)}_{k,i}(X)\mathbf{e}^{(j)}_k}$$
$$H_{i,k}^{(j)}=\alpha^{(j)}_{i,k}(X)\mathbf{e}^{(j)}_i+\alpha^{(j)}_{k,i}(X)\mathbf{e}^{(j)}_k$$
$$\lambda_{j,k}(X)=\displaystyle{\sum_{1\leq i \leq r_j,i\neq k}}
H_{i,k}^{(j)}.$$

\noindent As before, $X_N(t)=(V_N(t),\mathbf{e}_N(t))\in\mathbf{R}^p\times
\mathbf{R^d}$ is the solution of the stochastic hybrid model $(S_N)$.

\noindent Let $R_N(t)=\{(R_N^{(j)})_k(t)\}_{1\leq j\leq q,\ 1\leq k\leq r_j}$
with $R_N^{(j)}\in \mathbf{R}^{r_j}$ be defined as
$$(R_N^{(j)})_k(t)=\sqrt{N}\left(\{\mathbf{e}^{(j)}_N\}_k(t)-\{\mathbf{e}^{(j)}_N\}_k(0)-\int_0^t b_{j,k}(X_N(s))ds\right)$$

\begin{thm}

Under the same hypotheses as in \textbf{Theorem 2.1}, the process
$R_N$ converges in law, as $N\to \infty$, to the process
$R=\{(R^{(j)})_k(t)\}_{1\leq j\leq q,\ 1\leq k\leq r_j}$ with:
$$R^{(j)}(t)=\int_0^t \sigma^{(j)}(X(s))dW^{j}_s$$
where
\begin{itemize}
\item $X=(v,\mathbf{g})$ is the solution of the deterministic model $(D)$ with initial condition $X^{init}=X_N^{init}=X_0$,
\item $W^{j}$ are independent standard ${r_j}$-dimensional Brownian motions,
\item $\sigma^{(j)}(X)$ is the square root of matrix $G^{(j)}(X)$ s.t., for $1\leq k,l\leq r_j$:
$$\left\{
    \begin{array}{ll}
         G^{(j)}_{k,k}(X) = \lambda_{j,k}(X)\\
     G^{(j)}_{k,l}(X) = H^{(j)}_{k,l}(X) = G^{(j)}_{l,k}(X),\  l\neq k
    \end{array}
\right.$$
\end{itemize}
\end{thm}

\noindent This theorem leads to the following degenerate diffusion approximation $\tilde{X}_N=(\tilde{V}_N,\tilde{\mathbf{g}}_N) \in \mathbf{R}^p\times\mathbf{R}^d$, for $N$ sufficiently large:

$$(2.2)\  \  \left\{
    \begin{array}{ll}
         d\tilde{V}_N=f(\tilde{X}_N(t))dt\\
     d\mathbf{\tilde{g}}^{(j)}_{N}=b_{j}(\tilde{X}_N(t))dt+\frac{1}{\sqrt{N}}\sigma^{(j)}(\tilde{X}_N(t))dW^{j}_t
    \end{array}
\right.$$ where $b_j(X)$ is the vector $(b_{j,k}(X))_{1\leq k\leq
r_j}\in \mathbf{R}^{r_j}$, $1\leq j\leq q$.

\noindent Note that this approximation may not have the same properties as the original
process, even in the limit $N\to \infty$ (when considering for instance large
deviations \cite{note}).

\subsection{Functional central limit theorem and exit problem}

Let $X_N=(V_N,\mathbf{e}^{(1)}_N,..,\mathbf{e}^{(q)}_N)$ be the
solution of the stochastic model $(S_N)$ and
$X=(v,\mathbf{g}^{(1)},..,\mathbf{g}^{(q)})$ the solution of the
deterministic system $(D)$ with identical initial condition
$X^{init}=X_N^{init}=X_0 \in \mathbf{R}^{p+d}$.

\noindent Consider the $(p+d)$-dimensional processes:
$$Z_N=\left( \begin{array}{c} Y_N \\\mathbf{P}^{(1)}_N \\...\\\mathbf{P}^{(q)}_N  \end{array} \right) := \sqrt{N} \left( \begin{array}{c} V_N-v\\ \mathbf{e}^{(1)}_N-\mathbf{g}^{(1)} \\...\\ \mathbf{e}^{(q)}_N-\mathbf{g}^{(q)} \end{array} \right)$$
\noindent With this setting, we have the following result:
\begin{thm}
Under the same hypotheses as in \textbf{Theorem 2.1} the process
$Z_N$ converges in law, as $N\to \infty$ to the process $$Z=\left(
\begin{array}{c} Y \\\mathbf{P}^{(1)} \\...\\\mathbf{P}^{(q)}
\end{array} \right)$$ whose characteristic function
$\Psi(t,\theta)=\mathbf{E}\left[e^{i<\theta,Z(t)>}\right]$ satisfies
the following equation:

$$\frac{\partial\Psi}{\partial t} =\sum_{j=1}^q \left\{\sum_{l\in L}\sum_{k=1}^{r_j}\theta^{(j)}_k\frac{\partial b_{j,k}}{\partial x_l}\frac{\partial\Psi}{\partial \theta_l} -\frac{1}{2}\sum_{k,l=1}^{r_j}\theta^{(j)}_{k}\theta^{(j)}_lG^{(j)}_{k,l}\Psi\right\}+\sum_{m=1}^p\sum_{l\in L}\theta_m\frac{\partial f^m}{\partial x_l}\frac{\partial\Psi}{\partial \theta_l}$$

where $G^{(j)}_{k,l}$, $\frac{\partial f^m}{\partial x_l}$ and
$\frac{\partial b_{j,k}}{\partial x_l}$ are evaluated at $X(t)$, for
$1\leq m\leq p,\ 1\leq j\leq q,\ 1\leq k\leq r_j$ and $l\in L$,
with $\theta=((\theta_m)_{1\leq m\leq p},(\theta^{(j)}_k)_{1\leq
j\leq q,\ 1\leq k\leq r_j})=(\theta_l)_{\l\in L}$, and $L=\{(m)_{1
\leq m\leq p},\ (j,k)_{1\leq j\leq q,\ 1\leq k\leq r_j}\}$.

\end{thm}

\noindent Let us define $\tilde{Z}_N:=\left(\tilde{Y}_N,\tilde{P}_N\right):=\sqrt{N}\left(\tilde{X}_N-X\right)\in \mathbf{R}^p\times \mathbf{R}^d$, where we recall that $\tilde{X}_N$ is the Langevin approximation defined in $(2.2)$. Thus, for $1\leq j \leq q$, and $1\leq k \leq r_j$:
\begin{eqnarray*}
d\tilde{Y}_N &=&\sqrt{N}(f(\tilde{X}_N)-f(X))dt\\
d\tilde{P}_N^{j,k}&=&\sqrt{N}(b_{j,k}(\tilde{X}_N)-b_{j,k}(X))dt+\sigma^{(j)}(\tilde{X}_N)dW^j_t
\end{eqnarray*}
As an asymptotic linearization of $\tilde{Z}_N$, we define the diffusion process $\Theta=\left(\gamma, \pi \right) \in \mathbf{R}^p\times \mathbf{R}^d$ by:
\begin{eqnarray*}
d\gamma_m&=&\displaystyle{\sum_{l\in L}}\frac{\partial f^m}{\partial x_l}\Theta_l dt\\
d\pi_{j,k}&=&\displaystyle{\sum_{l\in L}} \frac{\partial
b_{j,k}}{\partial
x_l}\Theta_ldt+\displaystyle{\sum_{k'=1}^{r_j}}\sigma^{(j)}_{k,k'}(X)dW^{j,k'}_t
\end{eqnarray*}
\begin{thm}
The processes $Z$ and $\Theta$ have the same law.
\end{thm}
\noindent A computation of the moments equations for the limit process $Z$ for the example of section 2.1 is provided in Appendix B.

\noindent Using the central limit theorem 2.3, we provide in the next theorem a characterization of the fluctuations of first-exit time and location for the stochastic hybrid process $X_N$. Let $\phi:\mathbf{R}^{p+d} \to \mathbf{R}^{p+d}$ be continuously differentiable. Define
$$\tau_{N}:=\inf\{t \geq 0; \phi(X_N(t)) \leq 0\}$$
$$\tau := \inf\{t \geq 0; \phi(X(t)) \leq 0\}$$
the first passage times through $\phi=0$ respectively for the stochastic hybrid process and for its deterministic limit.
\begin{thm}
Assume the initial condition $X(0)$ satisfies $\phi(X(0))>0$. Suppose $\tau<\infty$ and $\nabla \phi(X(\tau)).F(X(\tau)) <0$. Denote the random variable $$\pi(\tau):= -\frac{\nabla \phi(X(\tau)).Z(\tau)}{\nabla \phi(X(\tau)).F(X(\tau))}$$ Then the following convergences in law hold when $N\to \infty$:
$$\sqrt{N}(\tau_N-\tau) \to \pi(\tau)$$
$$\sqrt{N}(X_N(\tau_N)-X(\tau)) \to Z(\tau)+\pi(\tau)F(X(\tau))$$
\end{thm}

\section{Application to neuron models}
In this section, we show how our previous theorems can be applied to
standard neuron models taking into account ion channel stochasticity.
\subsection{Kinetic formalism in neuron modelling}
Kinetic models can be used in many parts of nervous system
modelling, such as in ion channel kinetics, synapse and
neurotransmitters release modelling. Deterministic kinetic equations
are obtained as a limit of discrete stochastic models (hybrid or
not) as the population size, often the number of channels, is large.

As already mentioned in the introduction, compared to our general
formalism, the stochastic individuals are the ion channels and the
global variable $V_N$ the voltage potential. Constituted of several
subunits called gates, voltage-gated ion channels are metastable
molecular devices that can open and close. There exist different types of channels 
according to the kind of ions, and they are distributed within the
neuron membrane (soma, axon hillock, nodes of Ranvier, dendritic
spines) with heterogeneous densities.

In what follows, we consider the model of Hodgkin and Huxley, which
has been extended in different ways to include stochastic ion
channels.  In numerical studies, different versions have been
used, from a
\textit{two-state gating} interpretation \textit{e.g.} \cite{skaugen1979fbs} to
a \textit{multistate Markov}
scheme \cite{defelice1993,rowat2007iis}. In \cite{Steinmetz2000}, two of these models are compared,
one with a complete multistate Markov scheme, and the other inspired
from \cite{mainen1995msi} with a multistate scheme for the sodium
ion and a two-state gating for the potassium ion. We are here
considering only single-compartment neuron, but in order to deal
with spatial heterogeneities of axonal or ion channels properties
for instance, multi-compartmental models can be introduced as a
discretized description of the spatial neuron, with Ohm's Law
coupling between compartments.

\subsection{Application of the law of large numbers to Hodgkin and Huxley model}
Classically, the Hodgkin-Huxley model (HH) is the set of non-linear differential equations (3.1-3.4) which was introduced in \cite{HH52} to explain the ionic mechanisms behind action potentials in the squid giant axon. 
\begin{eqnarray}
\ C_m\frac{dV}{dt}&=& I - g_L(V-V_L) - g_{Na}m^3h(V-V_{Na}) - g_Kn^4(V-V_K)\\
\frac{dm}{dt}&=& (1-m)\alpha_m(V)-m\beta_m(V)\\
\frac{dh}{dt}&=& (1-h)\alpha_h(V)-h\beta_h(V)\\
\frac{dn}{dt}&=& (1-n)\alpha_n(V)-n\beta_n(V)
\end{eqnarray}
where $I$ is the input current, $C_m=1\mu F/cm^2$ is the a capacitance corresponding to the lipid bilayer of the membrane, $g_L=0.3 mS/cm^2$,$g_{Na}=120mS/cm^2$, $g_K=36ms/cm^2$ are maximum conductances and $V_L=10.6mV$, $V_{Na}=115mV$, $V_K=-12mV$ are resting potentials, respectively for the leak, sodium and potassium currents. The functions $\alpha_x,\ \beta_x$ for $x=m,n,h$ are opening and closing rates for the voltage-gated ion channels (see \cite{HH52} for their expression). The dynamics of this dynamical system can be very complex as shown in $\cite{guckoliva}$, but for our purpose let us describe only schematically the behavior of this system as the parameter $I$ is varied (see \cite{refHH} for more details). First for all $I$ there exists a unique equilibrium point. For $0<I<I_1\approx 9.8 \mu A/cm^2$, this equilibrium point is stable, and for $I_0<I<I_1$ where $I_0\approx 6.3 \mu A /cm^2$ this equilibrium coexists with a stable limit cycle and possibly many unstable limit cycles. At $I=I_1$ and $I=I_2$ occur two Hopf bifurcations. For $I_1<I<I_2\approx 153 \mu A/cm^2$, the equilibrium point is unstable and coexists with a stable limit cycle. For $I>I_2$, there are no more limit cycles, and the equilibrium point is stable. This bifurcation structure can be roughly interpreted as follows : for $I$ small the system converges to an equilibrium point, and for $I$ sufficiently large, the system admits a large amplitude periodic solution, corresponding to an infinite sequence of action potentials or spikes and the spiking frequency is modulated by the input current $I$. 

\noindent There are two stochastic interpretations of the Hodgkin-Huxley model involving either a multistate Markov model or a two state gating model. We now present them briefly and we apply our theorems to each of them. A comparison of these deterministic limits obtained for these models
is given in Appendix A and establishes an equivalence between the
deterministic versions as soon as initial conditions satisfy a
combinatorial coincidence relationship. This question is further studied in \cite{keener}, where the reduction of the law of jump Markov processes to invariant manifolds is investigated.
\paragraph{Multistate Markov model} This model has two types of ion channels : one for sodium and the other for potassium. The kinetic scheme describing the Markov jump process for one single
potassium channel is the following:

\includegraphics[height=1.3cm]{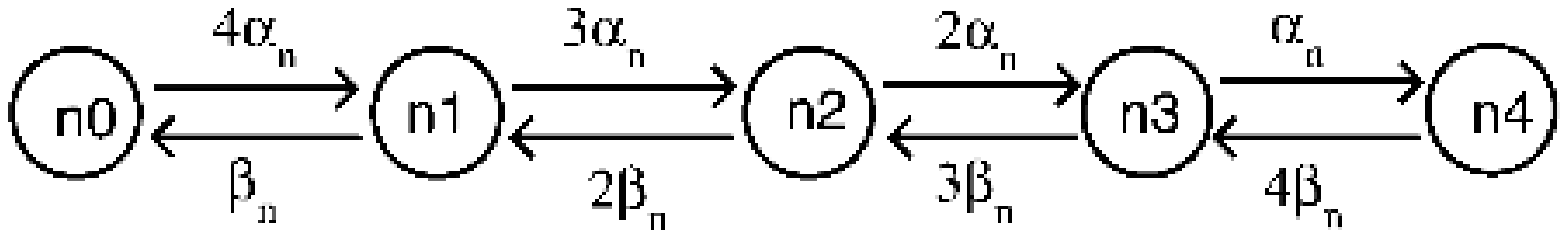}

\noindent And for sodium channel:

\includegraphics[height=3cm]{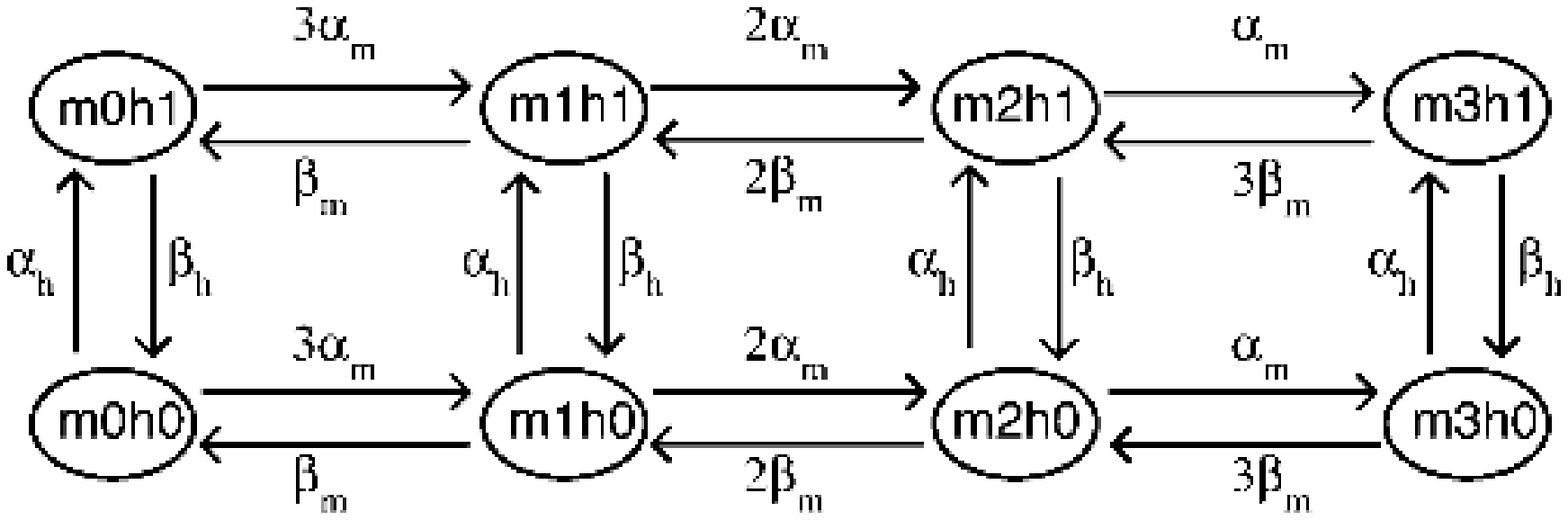}

\noindent All the coefficients in these two schemes are actually functions of
the membrane potential, and can be found in \cite{HH52}. The state spaces are: $$E_1=\{n_0,n_1,n_2,n_3,n_4\}$$ $$E_2=\{m_0h_1,m_1h_1,m_2h_1,m_3h_1,m_0h_0,m_1h_0,m_2h_0,m_3h_0\}$$
\noindent The
open states are respectively $n_4$ ($r_1=5$) and $m_3h_1$
($r_2=8$). The proportion of open potassium channels is denoted by
$u^{(1)}_N:=\{e_N^{(1)}\}_{n_4}$ and the proportion of open sodium channels by
$u^{(2)}_N:=\{e_N^{(2)}\}_{m_3h_1}$. In this model, The membrane potential dynamic is given by the
equation:$$\dot{V}_N(t)=-g_{Na}u^{(2)}_N(t)(V_N(t)-V_{Na})-g_{K}u^{(1)}_N(t)(V_N(t)-V_{K})-g_L(V(t)-V_L)+I$$
where $I\in \mathbf{R}$ is a constant applied current. With the notations of the previous sections, $f(v,u^{(1)}, u^{(2)})=-g_{Na}u^{(2)}(v-V_{Na})-g_{K}u^{(1)}(V_N(t)-V_{K})-g_L(v-V_L)+I$ and for example, $\alpha_{k,j}^{(1)}(v)=3\alpha_n(v)$ if $k=n_1,\ j=n_2$ and $\alpha_{k,j}^{(2)}(v)=2\beta_m(v)$ if $k=m_2h_0,\ j=m_1h_0$.

\noindent Applying \textbf{Theorem 2.1},  we obtain a
deterministic version of the stochastic Hodgkin-Huxley model when
$N\to\infty$, provided we choose the initial conditions appropriately:
$$\left\{
    \begin{array}{llll}
         \dot{v}=-g_{Na}e^{(2)}(t)(v(t)-V_{Na})-g_{K}e^{(1)}(t)(v(t)-V_K)-g_L(V(t)-V_L)+I\\
      \dot{\mathbf{g}}^{(j)}_k=\displaystyle{\sum_{1\leq i \leq r_j,i\neq k} \alpha^{(j)}_{i,k}(v)\mathbf{g}^{(j)}_i-\alpha^{(j)}_{k,i}(v)\mathbf{g}^{(j)}_k},\  \forall 1 \leq j \leq 2,\ \forall 1\leq k \leq r_j\\
     V(0)=v_{init}\\
     \mathbf{g}^{(j)}(0)=\mathbf{g}^{(j)}_{init}
    \end{array}
\right.$$
\noindent with $e^{(j)}=\mathbf{g}_{r_j}^{(j)}$, and where the rate functions
$\alpha^{(j)}_{m,p}$ are given in the above schemes, for
$j\in\{1,2\}$.

\paragraph{Two-state gating model} Another way of building a stochastic Hodgkin-Huxley model is to consider that the channels can be decomposed into independent gates. Each gate can be either open (state 1) or closed (state 0):
\begin{figure*}[!h]
\center
\includegraphics[scale=5]{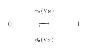}
\end{figure*}

with $z\in\{m,n,h\}$. The channel is open when all gates are open.\\
So, here, $q=3$ and $E_1=E_2=E_3=\{0,1\}$. If we denote
$u^{(z)}_N(t):=\{e_N^{(z)}\}_1$ the proportion of open gates $z$, for
$z\in\{m,n,h\}$, the membrane potential dynamic is then given
by:
\begin{eqnarray*}
\dot{V}_N(t)&=&-g_{Na}(u^{(m)}_N(t))^3u^{(h)}_N(t)(V_N(t)-V_{Na})\\
&-& g_{K}(u^{(n)}_N(t))^4(V_N(t)-V_{Na})-g_L(V_N(t)-V_L)+I
\end{eqnarray*}
\noindent which corresponds to $f(v,u^{(m)}, u^{(h)}, u^{(n)})=-g_{Na}(u^{(m)})^3u^{(h)}(v-V_{Na}) -g_{K}(u^{(n)})^4(v-V_{Na})-g_L(v-V_L)+I$.
In Figure 1, we give a sample trajectory of this two-state gating
stochastic Hodgkin-Huxley system. 

\noindent Applying \textbf{Theorem 2.1} gives the classical formulation of the
$4$-dimensional Hodgkin-Huxley model:
$$\left\{
    \begin{array}{lll}
         \dot{v}=-g_{Na}(u^{(m)}(t))^3u^{(h)}(t)(v(t)-V_{Na})-g_{K}(u^{(n)}(t))^4(v(t)-V_{Na})\\-g_L(V(t)-V_L)+I\\
     \dot{u}^{(z)}(t)=(1-u^{(z)}(t))\alpha_z(v(t))-u^{(z)}(t)\beta_z(v(t)),\  z\in\{m,n,h\}
    \end{array}
\right.$$

\begin{figure}[!h]
\centering
\includegraphics[height=11cm]{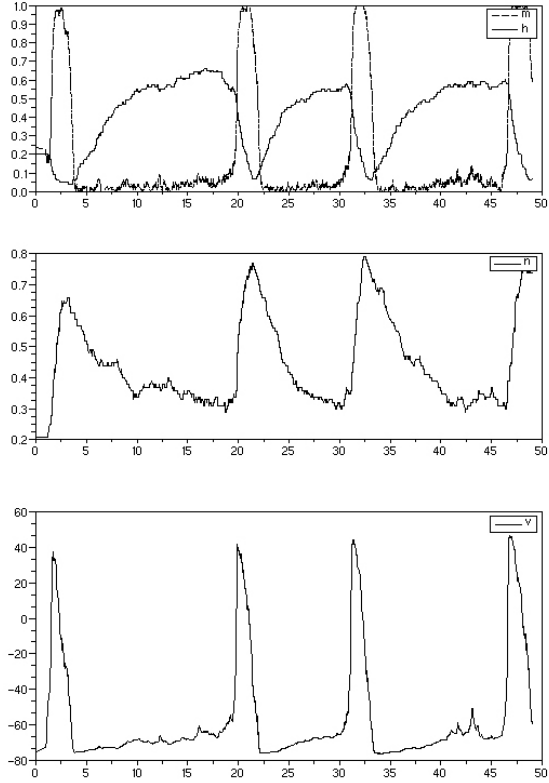}
\caption{Illustrative sample trajectory of a two-state gating
stochastic Hodgkin-Huxley model with $N=20$ (cf. section 3.). Top:
variables $m,h$ for the sodium channel (without unit). Middle:
variable $n$ for potassium channel (without unit). Down: variable
$V$ for membrane potential (unit: mV). Abscissa: time (arbitrary units).
Since $m,n$ and $h$ correspond to proportions of open gates, if one of
them is equal to 1, it means that all the corresponding gates are
open. An increase in the membrane potential $V$ causes an increase
in the proportion of open $m$ (sodium) gates, which in turn implies
an increase of $V$. This positive feedback results in a spike
initiation. Meanwhile, a further increase of $V$ is followed by a
decrease of the deactivation variable $h$, which closes the sodium
channels. This inhibition effect acts at a slower time-scale,
enabling a decrease of $V$. This decrease is strengthened by the
dynamic of variable $n$ (proportion of open potassium gates).}
\end{figure}

\subsection{Exponential convergence speed} We illustrate by numerical simulations the upper bound obtained in $(2.1)$ for the stochastic Hodgkin-Huxley
model with a two-state gating scheme. The number of sodium channels
$N_{Na}$ and potassium channels $N_K$ are proportional to the area
$S$ of the membrane patch. Thus, instead of $N$, $S$ will denote the
size parameter. For the squid giant axon, the estimated densities
for the ion channels used in the simulations are $\rho_{Na}=60
\mu m^{-2}$ and $\rho_{N_K}=18 \mu m^{-2}$.

\noindent We now display the results of numerical simulations of
$$P_S(T,\Delta)=\mathbf{P}\left[\displaystyle{\sup_{0\leq t \leq T} ||V_S(t)-v(t)||^2+ \sum_{j=1}^q ||\mathbf{e}^{(j)}_S(t)-\mathbf{g}^{(j)}(t)||^2 } > \Delta\right]$$
\noindent using Monte-Carlo simulations. We recall that from $(2.1)$:
$$\displaystyle \limsup_{S \to \infty} \frac{1}{S}\log P_S(T,\Delta) \leq - \frac{\Delta e^{-BT^2}}{CT} = C(T)$$
\noindent In Fig. 2, the simulation estimations of
$C_S(T)=\frac{1}{S}\log P_S(T,\Delta) $ are shown for different
values of $T$ and $S$  and can be compared to the
theoretical bound $C(T)$. Simulations are made without
input current, meaning that the stochastic solution is supposed to
fluctuate around the equilibrium point of the deterministic system in a
neighborhood of size proportional to $S^{-1/2}$. When $S$ increases, the simulation curve $C_S(T)$ is expected to pass below
the theoretical bound $C(T)$.
\begin{figure}[!h]
\centering
\includegraphics[scale=1]{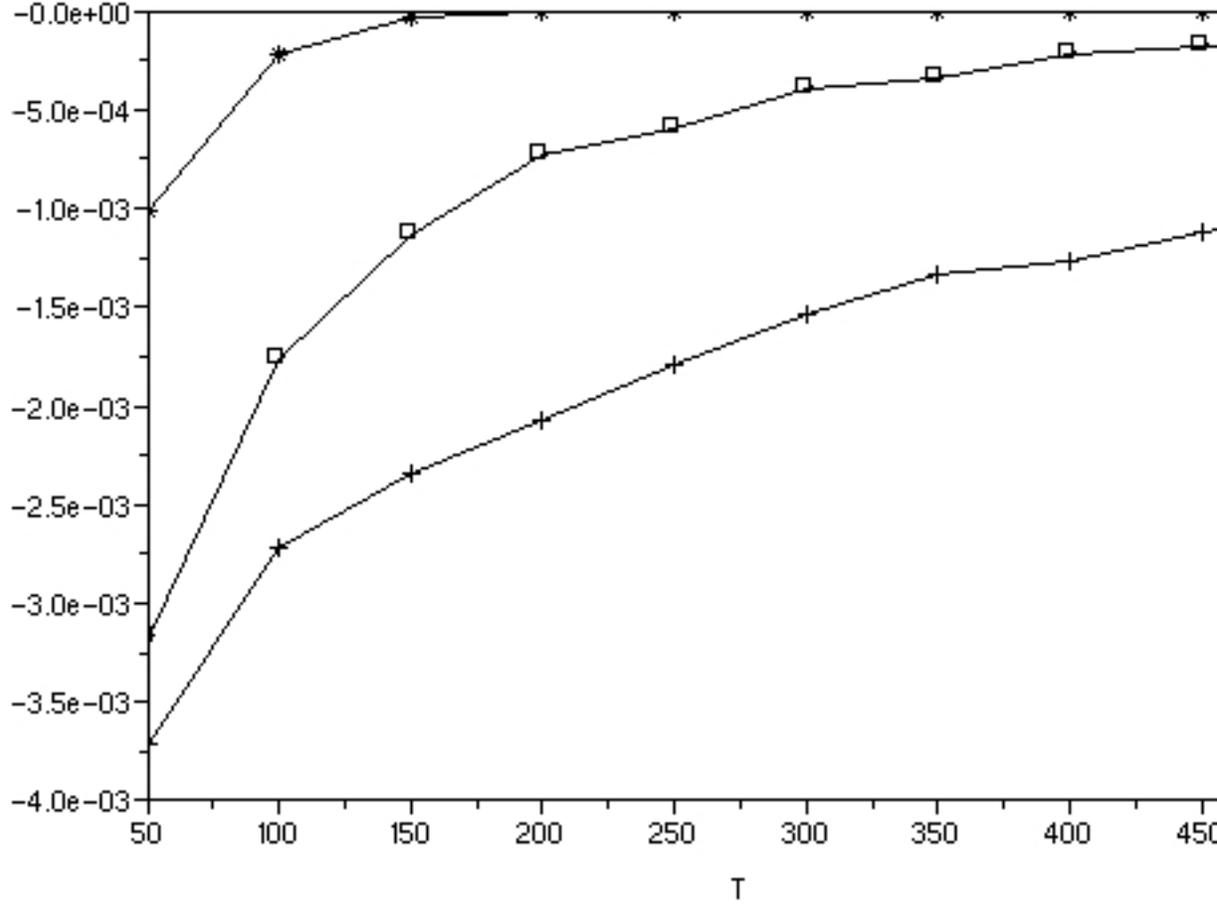}
\caption{Simulation results of the
Hodgkin-Huxley model with a ``two-state gating'' scheme with input
current $I=0$: this figure shows the quantity $\frac{1}{S}\ln
P_S(T,\Delta)$ as a function of $T$, where $S$ is the area of the
patch, and is thus proportional to $N$. Stars : $S=250\mu m^2$
(corresponding to $N_{Na}=15000$ and $N_K=5000$) ; Empty boxes :
$S=500\mu m^2$; Crosses: $S=750\mu m^2$. Lines are guide for the
eye.}
\end{figure}
\noindent For higher input currents, still subthreshold
($I<I_c$), but close to the bifurcation, channel noise will induce
spontaneous action potentials. For appropriate $\Delta$, the
probability $P_S(T,\Delta)$ can be interpreted as the probability
that the first spontaneous action potential (SAP) occurs before time
$T$. Thus the convergence speed bound gives an upper bound of the
repartition function of this first SAP time.

\noindent For higher input currents $I>I_c$, the
deterministic solution will be attracted by a stable limit cycle,
which corresponds to repetitive action potentials. In this case,
channel noise can introduce a jitter in the spiking times. Thus, if
one considers the supremum of the errors between the stochastic and
the deterministic solutions, this supremum will be quite large
(approximately the size of an action potential) as soon as the
difference between the stochastic spiking times and the
deterministic ones is of order the time course of an action
potential (2 ms). Thus, the supremum of the difference is not appropriate here and we will see in the following section how to quantify the impact of channel noise on the spiking frequency.

\subsection{Application of the central limit theorems}

In this section, we show how to investigate the
fluctuations around a stable fixed point (sub-threshold
fluctuations) and the fluctuations around a stable limit cycle
(firing rate fluctuations) using Theorem 2.3.
Let us consider a class of two-dimensional models,
corresponding to the \textbf{Example} of section 2.1. This class
contains reductions of the previous two-state gating Hodgkin-Huxley
model, or other models such as the Morris-Lecar model
\cite{morris1981}. Consider the process, with the notations of the \textbf{Example}:
$$\left( \begin{array}{c} Y_N \\ P_N \end{array} \right) = \left( \begin{array}{c} \sqrt{N}\left(V_N-V\right) \\ \sqrt{N}\left(e_N-g\right) \end{array} \right)$$
\noindent with initial conditions $(P_N(0),Y_N(0)) = (0,0)$. Then the 2-dimensional process $Z_N=\left(P_N,Y_N\right)$ converges
in law, as $N \to \infty$, towards the process $Z=\left(P,Y\right)$,
whose characteristic function is given by:
$$\mathbf{E}\left[e^{i(\theta_1P(t)+\theta_2Y(t))}\right]=e^{\theta_1^2A_t+\theta_2^2B_t+\theta_1\theta_2C_t}$$
\noindent Thus, defining $\Sigma_s$ the square root matrix of $\Gamma_s:=\begin{pmatrix} A'_s & C'_s/2 \\ C'_s/2 & B'_s \end{pmatrix}$, for
$0\leq s\leq t$, $Z$ can be written as a gaussian diffusion process:

$$Z_t=\int_0^t\Sigma_s dW_s$$ where $W$ is a standard two-dimensional brownian motion \footnote{The condition that the matrix $\Gamma_s$ admits a real square root
matrix can be reduced to $A'_s+B'_s \leq 0$ because one can show
that $\det(\Gamma_s)=A'_sB'_s-C'_s{}^2/4=0$ for all $s\geq 0$. This
condition is thus always satisfied because : $A'_0+B'_0 \leq 0$,
$A'_s$ and $B'_s$ have the same sign, and $(A'_s,B'_s,C'_s)$ cannot
cross $(0,0,0)$ by uniqueness of the solution of $z'=Mz$ (satisfied
by $y'$). The computation of the matrix $\Sigma_s$ gives:
$$\Sigma_s=\frac{\sqrt{-2(A'_s+B'_s)}}{A'_s+B'_s}\Gamma_s $$}.

\noindent From the equation for the characteristic function obtained in
\textbf{Theorem 2.3}, one derives that the triple $y=(A,B,C)$ is
solution of the system $\dot{y}_t=M_t y_t+E_t$ defined as:
$$\left( \begin{array}{c} \dot{A}_t \\ \dot{B}_t \\ \dot{C}_t \end{array} \right)= \left( \begin{array}{ccc}2b'_u & 0 & b'_v \\ 0 & 2f'_v & f'_u \\ 2f'_u & 2b'_v & b'_u+f'_v \end{array} \right) \left( \begin{array}{c} A_t \\ B_t \\ C_t \end{array} \right) + \left( \begin{array}{c} -\frac{1}{2}\lambda(V,u) \\ 0 \\ 0 \end{array} \right)\  \  (\mathbf{M})$$

\noindent with initial conditions $\left(0,0,0\right)$, and $\lambda(v,u)=\sqrt{(1-u)\alpha(v)+u\beta(v)}$. The partial derivatives $f'_v,\ f'_u,\ b'_v,\ b'_u$ and $\lambda$ are
evaluated at the deterministic solution $(V_t, g_t)$ .

\noindent We remark that, if $J$ is the Jacobian matrix at the point
$(V_t,g_t)$, and if its spectrum is
$\mbox{sp}(J)=\{\lambda_1,\lambda_2\}$ then the spectrum of $M$ is
$\mbox{sp}(M)=\{2\lambda_1,2\lambda_2,\lambda_1+\lambda_2\}$. Two different situations can be considered: 
\begin{itemize}
\item Starting from a fixed point $(V_0,u_0)$ of the deterministic system, the matrix $M_t=M(V_t,u_t)$ and the vector $E_t=E(V_t,u_t)$ are constant. One can derive an explicit analytical solution diagonalizing the matrix $M$. The time evolution for the variance and covariance of the difference between the deterministic solution and the stochastic one then depends on the stability $(\lambda_1,\lambda_2)$ of the considered fixed point.
\item Around a stable limit cycle (periodic firing): $M_t$ and $E_t$ are $T$-periodic functions. Using suitable coordinates and following Floquet's theory (see \cite{chicone99}), stability would be given by the spectrum of the solver $R(T):(A_0,B_0,C_0)\to(A_T,B_T,C_T)$. As explained in \cite{siam}, even if the real parts of the eigenvalues of the jacobian matrix are strictly negative for all time, unstable solutions may exist. In section 3.5 we investigate numerically the fluctuations around a stable limit cycle for the Morris-Lecar system.
\end{itemize}

\noindent If we consider $$\left( \begin{array}{c} \tilde{Y}_N \\ \tilde{P}_N
\end{array} \right) = \left( \begin{array}{c}
\sqrt{N}\left(\tilde{V}_N-V\right) \\
\sqrt{N}\left(\tilde{u}_N-u\right) \end{array} \right)$$ where
$(\tilde{V}_N,\tilde{u}_N)$ is the Langevin approximation, then the
moments equations, written for the linearized version around the
deterministic solution, give the same matrix $\Gamma_s$ at the limit
$N\to \infty$. But for finite $N$ the linearized process is not
gaussian (see Appendix B). Thus, our mathematical result can be directly related to the
simulations results obtained in \cite{Steinmetz2000}: in this paper
simulations of two neuron models with a large number of stochastic
ion channels are made, and the fluctuations of the membrane
potential below threshold exhibit approximately gaussian
distributions, but only for a certain range of resting potentials.
For smaller resting potentials, the shape of the distribution
remained unclear as it was more difficult to compute. Our approach
shows that, at finite $N$, for any range of the resting potentials the distribution is non-gaussian, but
when $N\to \infty$, the distribution tends to a gaussian, which
corresponds to the approximate gaussian distribution observed in the
simulations of \cite{Steinmetz2000}.

\subsection{Quantifying the effect of channel noise on neural coding}
Neurons encode incoming signals into trains of stereotyped
pulses referred to as action potentials (APs). It is the mean
firing frequency, that is the number of APs within a given time
window,  and the timing of the APs that are the main conveyors of
information in nervous systems. Channel noise due to the seemingly
random fluctuations in the opening and closing times of
transmembranar ion channels induces jitter in the AP timing and
consequently in the mean firing frequency as well. We show in the next subsections how our results can be applied to quantify these phenomena. The impact of channel noise on frequency coding is investigated in sec 3.5.1 and on spike timing coding in section 3.5.2. We close this section by some remarks concerning non-markovian processes arising when considering synaptic transmission in sec.3.5.3.

\subsubsection{Numerical study of the variance of spiking rate for Morris-Lecar model}
In this subsection, applying Theorem 2.3 to the Morris-Lecar system, we investigate the impact of channel noise on the variance of the firing frequency. The Morris-Lecar system was introduced in \cite{morris1981} to account for various oscillating states in the barnacle giant muscle fiber. We denote by $X=(V,m,n)$ the solution of:
\begin{eqnarray}
C_m\frac{dV}{dt}&=&I-g_L(V-V_L) - g_{Ca}m(V-V_{Ca}) - g_Kn(V-V_K):=F_v(X)\\
\frac{dm}{dt}&=&\lambda_m(V)(M_{\infty}(V)-m):=F_m(X)\\
\frac{dn}{dt}&=&\lambda_n(V)(N_{\infty}(V)-n):=F_n(X)
\end{eqnarray}
whrer $\lambda_m(V)=\cosh((V-V_1)/2V_2) $, $\lambda_n(V)= \phi_n*\cosh((V-V_3)/2V_4)$, $M_{\infty}(V)= (1 + \tanh [(V - V_1)/V_2)])/2$ and $N_{\infty}(V)= (1 + \tanh [(V - V_3)/V_4)])/2$.
We introduce as in the previous sections a stochastic version $X_N$ of this model with stochastic ion
channels, replacing the differential equation for $m$ and $n$ by birth-and-death processes with voltage-dependent opening rates
$\alpha_m=\lambda_m M_{\infty}$, $\alpha_n=\lambda_n N_{\infty}$ and closing rates $\beta_n=\lambda_n (1-N_{\infty})$. According to the parameters of the model, the
deterministic system $(3.5-3.7)$ may have a stable limit cycle $x^{LC}$
for some values of $I \in [I_{min},I_{max}]$ (see \cite{morris1981}). This corresponds to a phenomenon of regular spiking, characterized
by its rate. Assuming that the time length of a spike is almost constant, we suggest a proxy for this spiking rate:
\begin{equation*}
r(T):=\frac{1}{T}\int_0^T \phi_{th}(x(s)) ds
\end{equation*}
where $\phi_{th}$ is a sigmoid threshold function. In a similar way, we define the stochastic spiking rate by:
\begin{equation*}
r_N(T):=\frac{1}{T}\int_0^T \phi_{th}(X_N(s)) ds
\end{equation*}
\noindent As a candidate for $\phi_{th}$, we choose
$\phi_{th}(V):=\frac{e^{c(V-V_{th})}}{1+e^{c(V-V_{th})}}$ where $c$
and $V_{th}$ are two parameters.

\noindent A consequence of the central limit
theorem for $X_N$ is the following weak convergence:
\begin{equation*}
\sqrt{N}\left[r_N(T)-r(T)\right] \Rightarrow
R(T)=\frac{1}{T}\int_0^T Z(s).\nabla\phi_{th}(x(s))ds
\end{equation*}
where $Z$ is the weak limit of $\sqrt{N}\left[X_N-x\right]$:
\begin{equation*}
Z(s)=\int_0^s \Sigma(u)dW_u
\end{equation*}
$R(T)$ is a gaussian random variable with zero mean. For simplicity
we consider the case where $\phi_{th}$ is only a function of the
membrane potential $V$. Then the variance of $R(T)$ is:
\begin{equation}
\sigma_R^2(T)=\mathbf{E}\left[R(T)^2\right]=\frac{2}{T^2}\int_0^T
\int_0^s S_v(s')\phi'_{th}(V(s'))ds' \phi'_{th}(V(s))ds
\end{equation}
where $S_v(s)=\Sigma_{1,1}(s)$ is the variance of
$\sqrt{N}(V_N(s)-V(s))$.

\begin{figure}[!h]
  \center
  \includegraphics[scale=1.5]{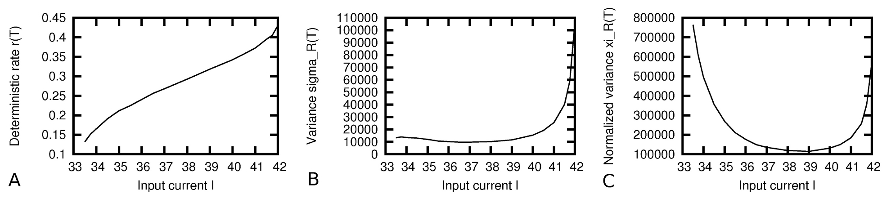}\\
  \includegraphics[scale=1.5]{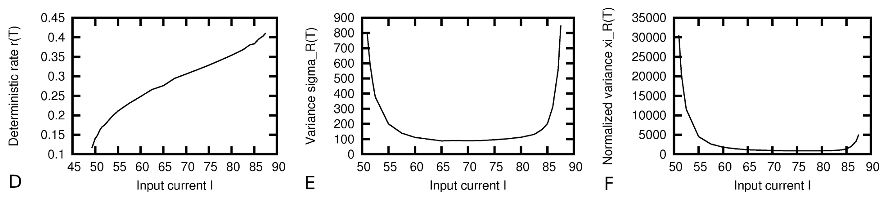}\\
  \caption{Impact of channel noise on the spiking rate.First row (ABC) : Class I regime. Second row (DEF) Class II regime. [A-D]. Deterministic rate $r(T)$ versus input current $I$ ($\mu A/cm^2$). [B-E]. Variance $\sigma_R^2(T)$ versus input current $I$ ($\mu A/cm^2$). [C-F]. Normalized variance $\xi(T)$ versus input current $I$ ($\mu A/cm^2$).  Parameters : for all figures $T=2000ms$, $c=10$, $V_{th}=0mV$, $C_m=20\mu F/cm^2$, $V_1=0mV$, $V_2=15mV$, $V_3=10mV$, $g_{Ca}=4mS/cm^2$, $g_K=8mS/cm^2$, $g_L=2mS/cm^2$, $V_K=-70mV$, $V_L=-50mV$, $V_{Ca}=100mV$, $\phi_n=0.1$ and $V_4=20mV$ for Class II (DEF) and $V_4=10mV$ for Class I (ABC).}
  \label{resultML}
\end{figure}

\noindent To estimate numerically the variance $\sigma_R^2(T)$, the first step is to determine numerically the limit cycle, then solve the moment equations (Appendix C) and immediately deduce $\Sigma(s)$. Thus the variance $\sigma_R^2$ can be computed using formula (3.8) without any stochastic simulation. In Fig. \ref{resultML} we show our numerical results, where we plot in C-F., as a
function of the input current $I$, the normalized variance $\xi(T)$ defined as:
\begin{equation*}
\xi(T) := \frac{\sigma_R^2(T)}{r(T)^2}
\end{equation*}

\paragraph{Comments} The value of $\xi(T)$ depends on a combination of the linear stability along the
cycle and on the variance of the noise (which is multiplicative)
along the cycle. If one wants to have the quantity $\frac{\mathbf{E}[(r_N(T)-r(T))^2]}{r(T)^2}$ of order $1$,
then the number $N$ of channels should be of order $\xi(T)$. Interestingly, this gives much smaller values for Class II than for Class I regime. In both cases, it corresponds to a reasonably small number of channels when $I$ is not too close from bifurcation points.

\subsubsection{Impact of channel noise on latency coding in the Morris-Lecar model}
\begin{figure}[!h]
  \center
  \includegraphics[width=12cm]{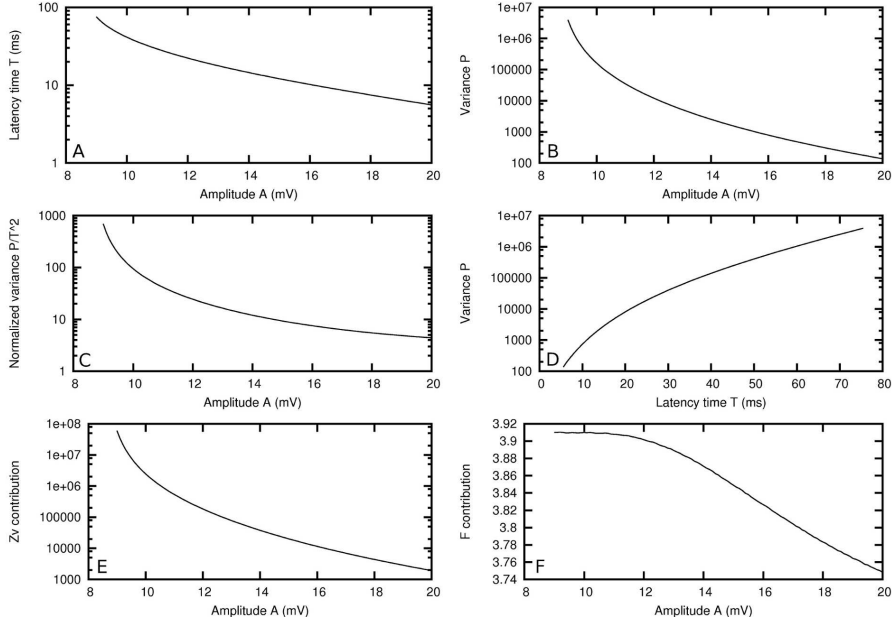}\\
  \caption{Impact of channel noise on latency coding. [A]. Latency time $T(A)$ versus amplitude $A$. [B]. Variance $P(A)$ versus amplitude $A$. [C]. Normalized variance $P(A)/T(A)^2$ versus amplitude $A$. [D]. Variance $P(A)$ versus latency time $T(A)$. [E]. Variance $S_v(A)$ versus amplitude $A$. [F]. Crossing speed $F_v(X(T(A)))$ versus amplitude $A$. Same parameters as in \ref{resultML} Class I, with input current $I=32\mu A/cm^2$.}
\label{latency}
\end{figure}

\noindent Whereas frequency coding requires an integration of the input signal over a relatively long time, individual spike time coding does not require such an integration. The time to first spike, called latency, depends on the value of the suprathreshold input. Thus it may have an interpretation in term of neural coding, and it has been shown in several sensory systems \cite{thorpe} that the first spike latency carries information. For example, a recent study \cite{gollish} concerning the visual system suggests that it allows the retina to transfer rapidly new spatial information. Impact of external noise on latency coding have been investigated in numerical studies \cite{pankratova} with stochastic simulations. We apply Theorem 2.5 to the Morris-Lecar model to investigate the impact of internal channel noise on first spike time. 
\noindent We chose the parameters (see \ref{resultML}) to obtain a Class I neuron model in the excitable regime. In this setting, there exists a unique steady state $X^*=(V^*, m^*, n^*)$. Starting from this equilibrium point, the impact of an input at $t=t_0$ is equivalent to an instantaneous shift of the membrane potential $V^* \to V^* + A$, where $A>0$ is the amplitude of this shift. Eventually the system goes back to its steady state, but if $A$ is higher than a threshold $A_{th}$ then a spike is emitted before going back to the steady state, whereas if $A$ is lower than $A_{th}$ no spike is emitted. For $A>A_{th}$, we define the latency time $T(A)$ as the elapsed time between $t_0$ and the spike. More precisely, let $(V_A(t),m_A(t),n_A(t))$ for $t\geq t_0$ be the solution of Morris-Lecar equations with initial conditions $X(t_0)=(V(t_0)=V^*+A, m(t_0)=m^*, n(t_0)=n^*)$. We define a spike as a passage of the membrane potential $V_A(t)$ through a threshold $V_{th}$. Then, with $t_0=0$ for simplicity, the latency time $T(A)$ can be written as $T(A):=\inf\{t\geq 0;\  V_A(t)>V_{th}\}$. As shown in Fig.\ref{latency}.A, the more $A>A_{th}$ is close to $A_{th}$, the longer is the latency time $T(A)$. The same setting can be extended in the stochastic case, defining a random variable $T_N(A)$. Applying Theorem 2.5, with $\phi(V,m,n)=V_{th}-V$, we express the variance $P(A)$ of the limit of $\sqrt{N}(T_N(A)-T(A))$ as $N\to \infty$:
\begin{equation}
P(A)=\frac{S_v(T(A))}{F_v(X(T(A)))^2}
\label{p}
\end{equation}
In (3.9), $S_v(T(A))$ is the variance of the $V$-component $Z_v$ of $Z$, where we recall that $Z$ is the limit of $\sqrt{N}(X_N-X)$ (see Theorem 2.3). The value of $S_v(T(A))$ is obtained from the numerical integration of the moments equations (..). The results are displayed in Fig.\ref{latency}, where the variance $P(A)$ and a normalized variance $P(A)/T(A)^2$ are plotted against the amplitude $A$ (\ref{latency}.B). In \ref{latency}.D the variance $P(A)$ is plotted against the latency time $T(A)$ (\ref{latency}.D). From (3.9), it appears that $P(A)$ is determined by two distinct contributions : the variance $S_v(T(A))$ (\ref{latency}.E) and the crossing speed $F(X(T(A)))$ (\ref{latency}.F) which actually does not influence much the variance $P(A)$. 
One way to interpret the results is the following: if $N$ is large, of order $P(A)$, then $\mathbf{E}[(T_N(A)-T(A))^2]$ is of order $1$. Thus, as an illustration, in order to keep $\mathbf{E}[(T_N(A)-T(A))^2]$ of order $1$, the required number of channels would be of order $10^2$ for a latency time of $10$ms and of order $10^5$ for latency time of $60$ms.

\subsubsection{Synaptic transmission and non-markovian processes}
In section 3.5.1, the quantity of interest was the firing frequency. However, the synaptic transmission between a neuron $A$ and a neuron $B$ has its own time scales. Therefore, neuron $B$'s input, called post-synaptic potential $\Psi^{A\to B}$, may be modeled as a functional of neuron $A$'s membrane potential $\{V_A(t)\}_{t\geq 0}$. Although synaptic transmission is presumably a non-linear process, one can consider as
a first approximation (cf. \cite{friesen}) that the process of interest is
obtained directly by the convolution of the process $V_A$ with some
kernel $K^{A\to B}$:

$$\Psi^{A\to B}(t)=\int_0^t K^{A\to B}(t-s) V(s)ds$$

\noindent The mathematical analysis of the impact of channel noise on this variable can be done in the light of theorems 2.1 and 2.3. Using the general notations for the stochastic process
and its deterministic limit, we define $\Psi_N(t) = \int_0^t K(t,s) X_N(s)ds$ and $\Psi(t) = \int_0^t K(t,s)x(s)ds$.
\paragraph{Law of large numbers}
Define $$S_N(T)=\displaystyle{\sup_{s\in[0,T]}}
|\Psi_N(t)-\Psi(t)|^2 $$ Clearly, using Cauchy-Schwartz inequality:
$$\mathbf{P}(S_N(T)>\Delta) \leq P_N(T,\eta(T)^{-1}\Delta)$$
with $$\eta(T)=T\displaystyle{\sup_{t\in[0,T]}} \int_0^t
|K(t,s)|^2ds$$ The convergence of $\Psi_N$ to $\Psi$ with the same
kind of exponential convergence speed is thus a direct consequence
of Theorem 2.1.

\paragraph{Gaussian fluctuations}
We know (Theorem 2.3) that $\sqrt{N}(X_N-x)$ converges weakly to the
diffusion $H(t)=\int_0^t R(u)dW_u$. As a consequence,
$\Omega_N=\sqrt{N}(\Psi_N-\Psi)$ converges also weakly, to the
following process:

$$\Omega(t)=\int_0^t K(t,s)\left(\int_0^s R(u)dW_u\right)ds$$

With an integration by part, one can rewrite:

$$\Omega(t)= \int_0^t Z(t,s)  dW_s$$
with
$$Z(t,s) = \int_s^t K(t,u)du R(s)$$

The process $\Omega$ is gaussian and one can easily compute
its variance as $\int_0^tZ(t,s)^2ds$. However, it is non markovian,
and some issues concerning the first hitting times of such processes
are solved in \cite{touboul}.

\section{Proof of the law of large numbers}

In this section we give the proof for \textbf{Theorem 2.1}. This
proof is inspired from \cite{Austin2007}, except for the
exponential martingale bound. In order to simplify the notation and
to make the arguments clearer and more intuitive, we write the proof
for the case of a single channel type with state space $\{0,1\}$ and
transition rates  given by the scheme:
\begin{figure*}[!h]
\center
\includegraphics[scale=5]{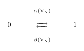}
\end{figure*}


In this case, the stochastic model $(S^0_N)$ is:
\begin{eqnarray*}
\dot{V}_N(t)=f(V_N(t),u_N(t)); V_N(0)=V_0\\
u_N(t)=\frac{1}{N}\sum_{k=1}^N \delta_1(u_t^{(k)}); u_N(0)=u^{(N)}_0
\end{eqnarray*}
where $u_t^{(k)}:0\to 1$ with rate $ \alpha(V_N(t))$ and $ 1\to 0$ with rate $\beta(V_N(t))$, for all $1\leq k\leq N$\\

The deterministic solution $(v,u)$ satisfies:
$$\left\{
    \begin{array}{lll}
        \dot{v}(t)=f(v(t),u(t))\\
        \dot{u}(t)=(1-u(t))\alpha(v(t))-u(t)\beta(v(t))\\
     v(0)=v_0\ ;\ u(0)=u_0
    \end{array}
\right.$$
In order to complete the proof, few slight changes in the notation
can be done:
\begin{itemize}

\item in order to work with more general jump Markov processes with finite state space, essentially all the expressions of the form $\delta_0(u)\alpha(v)-\delta_1(u)\alpha(v)$ should be replaced by
$$\sum_{i \neq j} \alpha_{i,j}(v,u)\delta_{e_i}(u)-\alpha_{j,i}(v,u)\delta_{e_j}(u)$$

\item in order to include $q$ different channel types (different ions), one should just write the same arguments for all the $q$ processes $\{\mathbf{e}^{(j)}_N(t)\}$ for $1\leq j \leq q$ and include all the $||\mathbf{e}^{(j)}_N(t)-\mathbf{e}^{(j)}(t)||$ for $1\leq j \leq q$ in the function $f(t)$ of Gronwall lemma application in section 3.4.

\end{itemize}

\subsection{Decomposition in a martingale part and a finite variation part}

\paragraph{Decomposition}
We decompose the difference between the stochastic and the deterministic
processes as a sum of a martingale part $M_N$ and a finite variation
part $Q_N$ as follows:
$$[u_N(t)-u_N(0)]-[u(t)-u(0)] = M_N(t)+\int_0^t Q_N(s)ds$$
where we define:
$$Q_N(t)=\frac{1}{N}\sum_{i=1}^N \left[\delta_0(u_t^{(i)})\alpha(V_N(t))-\delta_1(u_t^{(i)})\beta(V_N(t))\right]-\dot{u}(t)$$
$$M_N(t)= [u_N(t)-u(t)] -[u_N(0)-u(0)]-\int_0^t Q_N(s)ds $$

\paragraph{Lemma} As defined above, $(M_N(t))$ is a $\{F_t\}$-martingale.
\paragraph{Proof}
For $h>0$, define $\Delta M_N(t,h)=\frac{1}{h} \mathrm{E}
\left[M_N(t+h)-M_N(t) | F_t \right]$, then:
\begin{eqnarray*}
\Delta M_N(t,h) &=& \frac{1}{h}\frac{1}{N}\sum_{i=1}^N \mathrm{E}\left[\delta_1(u_{t+h}^{(i)})| F_t \right] -  \mathrm{E}\left[\delta_1(u_{t}^{(i)})| F_t \right]\\
&-& \frac{1}{h}\mathrm{E}\left[\int_t^{t+h}\left[\frac{1}{N}\sum_{i=1}^N \delta_0(u_s^{(i)})\alpha(V_s)-\delta_1(u_s^{(i)})\beta(V_s)\right]ds|F_t\right]\\
&-& \frac{1}{h}\left[u(t+h)-u(t)\right] + \frac{1}{h}\int_t^{t+h}
\dot{u}(s)ds
\end{eqnarray*}
\noindent The last line converges clearly to $0$ as $h \to 0$, and the two
first terms compensate as $h\to 0$. So we have:
$$\displaystyle{\lim_{h \to 0} \frac{1}{h} \mathrm{E} \left[M_N(t+h)-M_N(t) | F_t \right]} = 0$$
Therefore
$$\frac{d}{ds}\mathrm{E}[M_N(t+s)|F_t]|_{s=0}=0$$
By dominated convergence we have:
$$\frac{d}{ds}\mathrm{E}[M_N(t+s)|F_t]|_{s=s_0}=\mathrm{E}\left[\frac{d}{du}\mathrm{E}[M_{t+s_0+u}|F_{t+s_0}]|_{u=0}|F_t\right]=0$$
Finally:
$$\mathrm{E}[M_N(t+h)|F_t]=\mbox{Cst}=M_N(t)\  \  \   \otimes$$

\subsection{Martingale bound}
In this part we want to obtain a bound in probability for the
martingale part. We introduce the jump measure and the associated
compensator:

\noindent We define two random measures on $]0,T] \times \{0,1\}$:
\begin{itemize}
\item jump measure : $\kappa_i = \displaystyle{\sum_{t\in]0,T],u_t^{(i)}\neq u_{t^-}^{(i)}}} \delta_{(t,u_t^{(i)})}$
\item compensator :
$$\nu_i(dt,dy)=\left[\beta(V_N(t))\delta_1(u_{t^-}^{(i)})\delta_0(y)+\alpha(V_N(t))\delta_0(u_{t^-}^{(i)})\delta_1(y)\right]dt$$
\end{itemize}

We can rewrite $Q_N(s)$ and $M_N(t)$:
$$\int_0^t Q_N(s)ds=\frac{1}{N}\sum_{i=1}^N \int_{]0,T]\times\{0,1\}} (\delta_1(y)-\delta_1( u_{t^-}^{(i)}))\nu_i(ds,dy)-\int_0^t \dot{u}(s)ds$$
$$M_N(t)=\frac{1}{N}\sum_{i=1}^N \int_{]0,T]\times\{0,1\}} (\delta_1(y)-\delta_1( u_{t^-}^{(i)}))(\kappa_i-\nu_i)(ds,dy)$$

Then we have the following proposition:

\paragraph{Proposition}
Let $T>0$, $\epsilon>0$, $\delta >0$. Then there exists $N_0$ such
that $\forall N\geq N_0$,
$$\mathrm{P}\left[\displaystyle{\sup_{0\leq t \leq T} M_N(t)^2 } \geq \delta\right] \leq \epsilon$$
\paragraph{Proof}
Let us first recall that from standard results about residual processes (\cite{Kallenberg}) we have:
\begin{eqnarray*}
\mathrm{E}\left[M_N(t)^2\right]&=&\frac{1}{N^2}\sum_{i=1}^N \mathrm{E}\left[\int_{]0,T]x\{0,1\}} (\delta_1(y)-\delta_1( u_{t^-}^{(i)}))^2(\kappa_i-\nu_i)(ds,dy)\right]\\
&=&\frac{1}{N^2}\sum_{i=1}^N
\mathrm{E}\left[\int_{]0,T]}\beta(V_N(s))\delta_1(u_{s^-}^{(i)})+\alpha(V_N(s))\delta_0(u_{s^-}^{(i)})ds\right]
\end{eqnarray*}
Therefore, we can get a bound for
$\mathrm{E}\left[M_N(t)^2\right]$:
$$\mathrm{E}\left[M_N(t)^2\right] \leq C_1 \frac{t}{N} \max\left(||\alpha||_{\infty},||\beta||_{\infty}\right)$$
where $||\alpha||_{\infty}$ and $||\beta||_{\infty}$ are finite because $\alpha$ and $\beta$ are continuous and assumption (H1).
We then use Chebychev inequality and Doob inequality for $L^2$
martingales:
$$\mathrm{P}\left[\displaystyle{\sup_{0\leq t \leq T} M_N(t)^2 } \geq \delta\right] \leq \frac{1}{\delta} \mathrm{E}\left[\displaystyle{\sup_{0\leq t \leq T} M_N(t)^2}\right] \leq \frac{4}{\delta} \mathrm{E}\left[ M_N(t)^2\right]$$
and  $\mathrm{E}\left[ M_N(t)^2\right] \leq \frac{\epsilon \delta }{4}$ for all $N \geq N_0$.$\otimes$\\

In order to obtain a better estimate for the convergence rate, we
derive here an exponential bound for the martingale part. Our proof is inspired from techniques developed in \cite{darling2005}.

\paragraph{Proposition}
Let $T>0$,$\eta>0$. There exists a constant $C_{\eta}$ such that for
all $\delta \in ]0,\eta C_{\eta} T[$:
$$\mathrm{P}\left[\displaystyle{\sup_{0\leq t \leq T} |M_N(t)| } \geq \delta\right] \leq 2\exp\left(-\frac{\delta^2N}{2 C_{\eta}T}\right)$$

\paragraph{Proof}
We define, for $x=(u,v)$, $\theta \in R$:
\begin{eqnarray*}
m_N(x, \theta)&=&\int_R e^{\theta y} \lambda_N(x)\mu_N(x,dy) = N\lambda(x)[e^{\theta /N}\mu^{+}(x)+e^{-\theta /N}\mu^{-}(x)]
\end{eqnarray*}
\begin{eqnarray*}
\phi_N(x, \theta)&=&\int_R [e^{\theta y}-1-\theta y] \lambda_N(x)\mu_N(x,dy)\\
&=& \int_0^1 \frac{\partial^2m_N}{\partial
\theta^2}(x,r\theta)\theta^2(1-r)dr
\end{eqnarray*}
The second equality stems from integration by part. And, if $|\theta|<N \eta$,
\begin{eqnarray*}
|\frac{\partial^2m_N}{\partial \theta^2}(x,r\theta)|&=& \left|N\lambda(x)\frac{1}{N^2}[e^{r\theta /N}\mu^+(x)+e^{-r\theta /N}\mu^-(x)]\right|\leq \frac{C_{\eta}}{N}
\end{eqnarray*}
So, $|\phi_N(x, \theta)| \leq \frac{1}{2}\frac{C_{\eta}}{N}\theta^2$. Let us define
$$Z^{\epsilon}_N(t)=\exp[\epsilon M_N(t)-\int_0^t \phi_N((u_N(s),V_N(s)),\epsilon)ds]$$
$(Z^{\epsilon}_N(t))$ is a martingale thanks to
Dol\'eans Formula:
$$Z^{\epsilon}_N(t)=1+\int_0^t\int_R Z^{\epsilon}_N(s^-)[e^{\epsilon y}-1](\mu-\nu)(ds,dy)$$
Then we note $\tau=\inf\{t; M_N(t)>\delta \}$. On $\{\tau \leq t \}$,$Z^{\epsilon}_N(\tau) \geq \exp\{\delta \epsilon - \frac{t\epsilon^2 C_{\eta} }{2N}\}$. And by optional stopping theorem:
$$\mathrm{E}[Z^{\epsilon}_N(\min(t,\tau))]=\mathrm{E}[Z^{\epsilon}_N(0)] \geq \mathrm{E}[Z^{\epsilon}_N(\tau)\mathrm{1}_{\tau \leq t}] \geq \mathrm{P}(\tau \leq t)\exp\{\delta \epsilon - \frac{t\epsilon^2 C_{\eta}}{2N}\}$$
So, $\mathrm{P}\left[\displaystyle{\sup_{0\leq t \leq T} |M_N(t)| } > \delta \right] = \mathrm{P}(\tau \leq T) \leq \exp\{-\delta \epsilon + \frac{T\epsilon^2 C_{\eta} }{2N}\}$.\\Finally when $\delta \in ]0,\eta C_{\eta} T[$, with $\epsilon =
\frac{\delta N}{C_{\eta}t}$, and applying the same argument to $-M_N(t)$
we get the result.$\otimes$

\subsection{Finite Variation Part}
In this section we use the Lispchitz property of
$\alpha$ and $\beta$ to provide a bound for the finite
variation part, in order to apply later Gronwall Lemma.
\paragraph{Lemma}
There exists $C_1>0$ independent of $N$ such that:
$$|Q_N(t))| \leq C\left(|u_N(t)-u(t)|+|V_N(t)-v(t)|\right)$$
\paragraph{Proof}
\begin{eqnarray*}
Q_N(t) &=& \frac{1}{N}\sum_{i=1}^N\delta_0(u_t^{(i)})\alpha(V_N(t)) - (1-u(t))\alpha(v(t))\\
&-& \frac{1}{N}\sum_{i=1}^N\delta_1(u_t^{(i)})\beta(V_N(t)) -
u(t)\beta(v(t))
\end{eqnarray*}
Let us start with the second term of the difference, called $Q^{1\to
0}$:
\begin{eqnarray*}
Q^{1\to 0}&=& \frac{1}{N}\sum_{i=1}^N\delta_1(u_t^{(i)})\beta(V_N(t)) - u(t)\beta(v_t)\\
&=&\frac{1}{N}\sum_{i=1}^N\delta_1(u_t^{(i)})\beta(V_N(t)) - u(t)\beta(V_N(t))+ u(t)\left(\beta(V_N(t))-\beta(v(t))\right)\\
&=&\beta(V_N(t))\left( u_N(t)-u(t) \right) +
\underbrace{u(t)}_{\in[0,1]}\left(\beta(V_N(t)) - \beta(v(t))
\right)
\end{eqnarray*}
Then,
\begin{eqnarray*}
|Q^{1\to 0}|&\leq& ||\beta||_{\infty} |u_N(t)-u(t)| +
K_{\beta}|V-v(t)|
\end{eqnarray*}
where $K_{\beta}$ is the Lipschitz coefficient of $\beta$.
We do the same for the other term of the difference:
\begin{eqnarray*}
|Q^{0\to 1}|&\leq& ||\alpha||_{\infty} |u_N(t)-u(t)| +
K_{\alpha}|V_N(t)-v(t)|
\end{eqnarray*}
So the proof is complete, with $C_1=\max(||\alpha||_{\infty},
||\beta||_{\infty}, K_{\alpha}, K_{\beta})$ $\otimes$

If more general transition rates $\alpha(v,u)$ and $\beta(v,u)$ depend on $v$ \textit{and} $u$, one would need to replace $||\alpha||_{\infty}$ and  $||\beta||_{\infty}$, respectively by  $||\alpha||_{\infty}+K_{\alpha}^{(u)}$ and  $||\beta||_{\infty}+K_{\beta}^{(u)}$, where $K_{\alpha}^{(u)}, K_{\beta}^{(u)}$ are the Lipschitz coefficients associated with the second variable $u$.

\subsection{Proof of theorem 2.1}

\paragraph{Law of large numbers}
We want to apply Gronwall Lemma to the function:$$f(t)=
|V_N(t)-v(t)|^2+ |u_N(t)-u(t)|^2$$
From the previous section we have a good control on the martingale term and the following estimate:
\paragraph{Corollary}
There exists  $C_2>0$ independent of $N$ such that:
\begin{eqnarray*}
|u_N(t)-u(t)|^2&\leq&4[|u_N(0)-u(0)|^2 + C_2T\int_0^t|u_N(s)-u(s)|^2ds \\
&+& C_2T\int_0^t|V_N(s)-v(s)|^2+M_N(t)^2]
\end{eqnarray*}
\paragraph{Proof} As $u_N(t)-u(t)=u_N(0)-u(0) +M_N(t)+\int_0^tQ_N(s)ds$ and $(x+y+z)^2\leq 4(x^2+y^2+z^2)$, the result is a direct application of the previous lemma and of Cauchy-Schwarz inequality.$\otimes$
We need
now to work on $|V_N(t)-v(t)|^2$, using hypothesis (H1) ,
with $K_1=\displaystyle{\sup_N \sup_{s\in [0,T]}}
\left|\frac{\partial f}{\partial v}(V_N(s),u_N(s))\right|$ and
$K_2=\displaystyle{\sup_N \sup_{s\in [0,T]}} \left|\frac{\partial
f}{\partial u}(V_N(s),u_N(s))\right|$.

\noindent Between the jumps, we have:
\begin{eqnarray*}
\frac{d}{dt}\left( |V_N(t)-v(t)|^2\right) &=&
2\left(f(V_N(t),u_N(t))-f(v(t),u(t))\right)\left(V_N(t)-v(t)\right)
\end{eqnarray*}
Thus,
\begin{eqnarray*}
 |V_N(t)-v(t)|^2 &=& 2\int_0^t\left[f(V_N(s),u_N(s))-f(v(s),u(s))\right]\left(V_N(s)-v(s)\right)ds\\
&+& |V_N(0)-v_0|^2\\
&\leq& |V_N(0)-v_0|^2 + 2K_1\int_0^t |V_N(s)-v(s)|^2ds\\
&+& 2K_2\int_0^t |u_N(s)-u(s)||V_N(s)-v(s)|ds\\
&\leq& |V_N(0)-v_0|^2 + 2K_1\int_0^t |V_N(s)-v(s)|^2ds\\
&+& K_2\int_0^t |u_N(s)-u(s)|^2ds+K_2 \int_0^t|V_N(s)-v(s)|^2ds
\end{eqnarray*}
where we used successively Cauchy-Schwartz inequality and $ab \leq
\frac{1}{2}(a^2+b^2)$. Putting together this inequality with
the Corollary we obtain:
$$f(t) \leq A + B\int_0^tf(s)ds$$
where $B=B(T)=\max(2K_1(T)+K2(T),C_2T)$ does not depend on $N$ and
is linear w.r.t  $T$ if (H2) holds, and
$$A = |u_N(0)-u_0|^2 + |V_N(0)-v_0|^2 + K_A \displaystyle{\sup_{0 \leq s \leq T} M_s^2}$$
If we control the initial conditions, then, with the control we have
on the martingale part, $A$ can be chosen arbitrarily small (with
high probability) and we can conclude with Gronwall Lemma.

\paragraph{Exponential convergence speed}

If the initial conditions are the same for the
stochastic and deterministic model, we actually
have a exponentially fast convergence, thanks to the exponential bound for the martingale part:
there exists a constant
$C_m>0$ such that:
$$\displaystyle \limsup_{N \to \infty} \frac{1}{N}\log\mathrm{P}\left[\displaystyle{\sup_{0\leq t \leq T} |V_N(t)-v(t)|^2+ |u_N(t)-u(t)|^2 } > \Delta\right] \leq - \frac{\Delta e^{-B(T)T}}{2K_AC_mT}$$

\section{Proof of the central limit theorems}

As before, we write the proofs for the case of a single
channel type with state space $\{0,1\}$ and transition rates  given
by the scheme:
\begin{figure*}[!h]
\center
\includegraphics[scale=5]{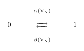}
\end{figure*}


\subsection{Langevin approximation}

In this case, Theorem 2.2 can be written as follows:

Let $b(u,v)=(1-u)\alpha(v)-u\beta(v)$, and $\left(V_N,u_N\right)$
solution of the stochastic model $(S_N)$. Then, the process $R_N(t)=\sqrt{N}\left(u_N(t)-u_N(0)-\int_0^t
b(u_N(s),V_N(s))ds\right)$ converges in law, as $N \to \infty$,
towards the process $R(t)$ defined as a stochastic integral:
$$R(t) = \int_0^t \sqrt{(1-u(s))\alpha(V(s)) + u(s)\beta(V(s))} dB_s$$
where $B$ is a standard brownian motion and $u(t),V(t)$ is the
unique solution of:
$$\dot{V}=f(V,u)$$
$$\dot{u}=(1-u)\alpha(V) - u\beta(V)$$
$$\forall N,\  u(0)=u_0=u_N(0)$$
$$\forall N,\  V(0)=V_0=V_N(0)$$

This result provides the following degenerate diffusion
approximation $(\tilde{V}_N,\tilde{u}_N)$, for $N$ sufficiently large:

$$d\tilde{V}_N(t)=f(\tilde{V}_N(t),\tilde{u}_N(t))dt$$
$$d\tilde{u}_N(t)=\left[(1-\tilde{u}_N(t))\alpha(\tilde{V}_N(t))-\tilde{u}_N(t)\beta(\tilde{V}_N(t))\right]dt+\sigma_N(\tilde{u}_N(t),\tilde{V}_N(t))dB_t$$
$$\sigma_N(u)^2=\frac{1}{N}\left[(1-u)\alpha(v)+u\beta(v)\right]=\frac{1}{N}\lambda(v,u)$$

\noindent Let
$g(u,v)=\lambda_N(u,v)[\frac{1}{N^2}\mu^+(u,v)+\frac{1}{N^2}\mu^-(u,v)]=(1-u)\alpha(v)+u\beta(v)$.

\noindent Note that in the multidimensional case, the real valuedfunction $g$ above
becomes a $d\times d$-matrix. Since the different
channel types $j$ are supposed to independent, this matrix would be bloc
diagonal, with blocs of size $r_j$, thus assuring the independence
of the $q$ ($r_j$-dimensional Brownian motions) $W^{(j)}$ in Theorem
2. The blocs of size $r_j$ are given by the matrix $G^{(j)}$ of
theorem 2, and arise from the calculation of the covariances:
$$G^{(j)}_{i,j}(x)=N\lambda_N(x)\int_E z_i z_j \mu_N(x,z)$$

\paragraph{Proof of Theorem 2.2}
We adapt the proof given by Kurtz \cite{Kurtz1971}: we prove the convergence of characteristic functions
plus tightness.
The tightness property follows from the inequality:
$$\mathrm{P}[\displaystyle{\sup_{s\leq T}|R_N(s)}|>\delta]\leq \frac{tN}{\delta^2}||g||_{\infty}$$
\noindent Let $\phi_N(t,\theta)=\mathbf{E}[e^{i\theta R_N(t)}]$ the characteristic function of $R_N$. Let $h(M_N(t))=e^{i\theta R_N(t)}$, $\sqrt{N}M_N(t)=R_N(t)$, $\psi(u)=\frac{e^{iu}-1-iu+u^2/2}{u^2}$, $\xi(u)=e^{iu}-1-iu=u^2\psi(u)-u^2/2$.
We then have:
\begin{eqnarray*}
\phi_N(t,\theta)-1 &=&\mathbf{E}[h( M_N(t))]-h(0)\\
&=& \int_0^t \mathbf{E}[\lambda_N(s)\int_{E_N}h(w-u_N(s)+M_N(s))-h(M_N(s))\\
&-&(w-u_N(s))h'(M_N(s))\mu_N(s,dw)]ds\\
&=&\int_0^t \mathbf{E}[e^{i\theta R_N(s)} \lambda_N(s)\int_{E_N}\xi(\theta\sqrt{N}(w-u_N(s))\mu_N(s,dw)]ds\\
&=&-\int_0^t \mathbf{E}\left[\frac{1}{2}e^{i\theta R_N(s)}\lambda_N(s)\int_{E_N} N\theta^2(w-u_N(s))^2\mu_N(s,dw)\right]ds\\
&+& \int_0^t \mathbf{E}[e^{i\theta R_N(s)}\lambda_N(s)\int_{E_N} N\theta^2(w-u_N(s))^2\\
&\times&\psi\left(\sqrt{N}\theta(w-u_N(s))\right)\mu_N(s,dw)]ds
\end{eqnarray*}
where $\lambda_N(s)$ stands for $\lambda_N(u_N(s),V_N(s))$ and
$\mu_N(s,dw)$ for $\mu_N(u_N(s),V_N(s),dw)$. The second term in the last equality,
call it $K_N(\theta)$, converges to $0$ as $N\to{\infty}$ by
dominated convergence, and because
$\psi\left(\sqrt{N}\theta(w-u_N(s))\right)=\psi\left(^+_-\theta /
\sqrt{N}\right) \to 0$ as $\displaystyle{\lim_{u\to 0}} \psi(u)=0$.
So we have:
\begin{eqnarray*}
\phi_N(t,\theta)-1 &=& -\int_0^t \mathbf{E}\left[\frac{1}{2}e^{i\theta R_N(s)}\theta^2g(u_N(s),V_N(s))\right]ds + K_N(\theta)\\
&=&-\frac{1}{2}\int_0^t\theta^2 g(u(s),V(s))\phi_N(s,\theta)ds\\
&+&\frac{1}{2}\int_0^t\theta^2\mathbf{E}\left[(g(u(s),V(s))-g(u_N(s),V_N(s)))e^{i\theta R_N(s)}\right]ds\\
&+& K_N(\theta)
\end{eqnarray*}

\noindent Again, the second term in the last equality, call it $J_N(\theta)$, converge to $0$ as $N\to \infty$, because of the convergence of $u_N$ and $V_N$ to $u$ and $V$.(cf. Theorem 2.1)\\
By Gronwall lemma, we  conclude that $\phi_N(t,\theta) \to
\phi(t,\theta)$ with:

$$\phi(t,\theta)=\exp\{-\frac{1}{2}\theta^2 \int_0^t g(u(s),V(s))ds\} \ \ \otimes$$

\subsection{Functional central limit theorem}

Let $(V_N,u_N)$ be the solution of the simplified stochastic model
$(S_N)$ and $(V,u)$ of the deterministic model $(D)$ introduced in the Example of section 2.
\noindent Consider the process:
$$\left( \begin{array}{c} P_N \\ Y_N \end{array} \right) = \left( \begin{array}{c} \sqrt{N}\left(u_N-u\right) \\ \sqrt{N}\left(V_N-V\right) \end{array} \right)$$
\noindent If the initial conditions satisfy $(P_N(0),Y_N(0)) = (0,0)$, the 2-dimensional process $\left(P_N,Y_N\right)$ converges in
law, as $N \to \infty$, towards the process $\left(P,Y\right)$,
with characteristic function:
$$\mathbf{E}\left[e^{i(\theta_1P(t)+\theta_2Y(t))}\right]=e^{\theta_1^2A(t)+\theta_2^2B(t)+\theta_1\theta_2C(t)}$$
\noindent The functions $A,B$ and $C$ are solutions of the system:
$$\left( \begin{array}{c} A' \\ B' \\ C' \end{array} \right)= \left( \begin{array}{ccc}2b'_u & 0 & b'_v \\ 0 & 2f'_v & f'_u \\ 2f'_u & 2b'_v & b'_u+f'_v \end{array} \right) \left( \begin{array}{c} A \\ B \\ C \end{array} \right) + \left( \begin{array}{c} -\frac{1}{2}\lambda(V,u) \\ 0 \\ 0 \end{array} \right)\  \ (\mathbf{M})$$
\noindent with initial conditions $\left(0,0,0\right)$, and with
$\lambda(v,u)=\sqrt{(1-u)\alpha(v)+u\beta(v)}$.

\paragraph{Proof of Theorem 2.3}

Just as in the proof of Theorem 2.2, let us define:
\begin{eqnarray*}
\phi_N(t,\theta)&=&\mathbf{E}\left[e^{i\left(\theta_1P_N(t)+\theta_2Y_N(t)\right)}\right]
\end{eqnarray*}
Let us also define $Z_N=(u_N-u,V_N-V)$, $X_N=(u_N,V_N)$, $X=(u,V)$, and
$h(x,y)=e^{i\sqrt{N}\left(\theta_1 x + \theta_2 y\right)}$. Then:
\begin{eqnarray*}
\phi(t,\theta)-1&=&\mathbf{E}\left[h(Z_N(t))-h(Z_N(0))\right]\\
&=& \int_0^t\mathbf{E}[N\lambda(X_N(s))\int_{E_N}\{h\left(w-u(s),V_N(s)-V(s)\right)\\
&-& h\left(Z_N(s)\right)\}\mu(X_N(s),dw)\\
&-& b(X(s))h'_x\left(Z_N(s)\right)+ \left(
f\left(X_N(s)\right)-f\left(X(s)\right)\right)h'_y\left(Z_N(s)\right)]ds
\end{eqnarray*}
So $\phi_N(t,\theta)-1  =  G_N(\theta,t)+H_N(\theta,t)$
with
\begin{eqnarray*}
G_N(\theta,t)&=&\int_0^t\mathbf{E}\left[\Omega_N(s)\left\{(e^{i\theta_1\sqrt{N}/N}-1)\mu_{+}+(e^{-i\theta_1\sqrt{N}/N}-1)\mu_{-}\right\}\right]ds\\
\Omega_N(s)&=&N\lambda(X_N(s)) h(Z_N(s)),\  
\mu_{+/-}=\mu_{+/-}(X_N(s))\\
H_N(\theta,t)&=&\int_0^ti\mathbf{E}\left[-\theta_1\sqrt{N}b(X(s))h(Z_N(s))+\theta_2\sqrt{N}\left\{f\left(X_N(s)\right)-f\left(X(s)\right)\right\}\right]ds
\end{eqnarray*}
Then in order to use the asymptotic development of $e^x$ when $x\to
0$ we introduce the function $K(u)=e^{iu}-1-iu+u^2/2$. Then,
knowing that $\mu_{+}+\mu_{-}=1$:
\begin{eqnarray*}G_N(\theta,t)=\int_0^t\mathbf{E}\left[\Omega_N(s)\left\{i\frac{\theta_1}{\sqrt{N}}(\mu_{+}-\mu_{-})(X_N(s))-\frac{\theta_1^2}{2N}+K(\theta_1/\sqrt{N})\right\}\right]ds
\end{eqnarray*}
\noindent Since $b(x)=\lambda(x)(\mu_{+}(x)-\mu_{-}(x))$, we have :
\begin{eqnarray*}
G_N(\theta,t)&=&\int_0^t\mathbf{E}\left[i\theta_1\sqrt{N}b(X_N(s))h(Z_N(s))\right]ds \\
& + & \int_0^t\mathbf{E}\left[-\frac{\theta_1^2}{2}\lambda(X_N(s))h(Z_N(s))\right]ds \\
& + & \int_0^t\mathbf{E}\left[NK(\theta_1/\sqrt{N})h(Z_N(s)\right]ds
\end{eqnarray*}
Therefore:
\begin{eqnarray*}
\phi_N(t,\theta)-1&=&\int_0^t \mathbf{E}\left[-\frac{1}{2}\theta_2^2\lambda(X_N(s))h(Z_N(s))\right]ds  \  (A)\\
&+&\int_0^t \mathbf{E}\left[h(Z_N(s))i\theta_1\sqrt{N}\left\{b(X_N(s))-b(X(s))\right\}\right]ds  \  (B)\\
&+&\int_0^t \mathbf{E}\left[h(Z_N(s))i\theta_2\sqrt{N}\left\{f(X_N(s))-f(X(s))\right\}\right]ds  \  (C)\\
&+&\int_0^t \mathbf{E}\left[h(Z_N(s))NK(\theta_1/\sqrt{N})\lambda(X_N(s))\right]ds  \  (D)
\end{eqnarray*}
Using the derivatives of $b$ and $f$, and the convergence of $X_N$
to $X$ we can make a development of the sum
$B+C$:
\begin{eqnarray*}
B+C &=& \int_0^t \mathbf{E}[h(Z_N)i\sqrt{N}\{(u_N-u)(\theta_1b'_u+\theta_2 f'_u)\\
&+& (V_N-V)(\theta_1b'_v+\theta_2 f'_v)\}]ds\}+\epsilon_N(t,\theta)
\end{eqnarray*}
where we dropped the \textit{s} and where $b'_u,b'_v, f'_u, f'_v$ are taken at $X_N(s)$.\\
Noting that $h(Z_N)i\sqrt{N}(u_N-u)=h'_x(Z_N)$ and
$h(Z_N)i\sqrt{N}(V_N-V)=h'_y(Z_N)$, we have:
$$B+C  =  \int_0^t \mathbf{E}\left[h'_x(Z_N)(\theta_1b'_u+\theta_2 f'_u)+h'_y(Z_N)(\theta_1b'_v+\theta_2 f'_v)\right]$$
And the term $D$ converges to zero as $N \to \infty$ by dominated convergence since $K(u)/u^2$ is bounded and converges to $0$.\\
As we have the convergence in Theorem 2.1 of $X_N$ to $X$,we get
the convergence of $\phi_N(t,\theta)$ to $\Psi(t,\theta)$,
satisfying:
\begin{eqnarray*}
\frac{\partial\Psi}{\partial t}(t,\theta) &=& -\frac{1}{2}\theta_1^2\lambda(X(t))\Psi(t,\theta)+(\theta_1b'_u(X(t))+\theta_2 f'_u(X(t)))\frac{\partial\Psi}{\partial \theta_1}\\
&+& (\theta_1b'_v(X(t))+\theta_2
f'_v(X(t)))\frac{\partial\Psi}{\partial \theta_2}
\end{eqnarray*}

Tightness stems from the Markov property and the following
estimate obtained in the proof Theorem 2.1:
\begin{eqnarray*}
\displaystyle\mathbf{P}\left[\displaystyle{\sup_{0\leq t \leq T} ||\sqrt{N}\left[V_N(t)-v(t)\right]||^2+ \sum_{j=1}^q ||\sqrt{N}\left[\mathbf{e}^{(j)}_N(t)-\mathbf{g}^{(j)}(t)\right]||^2 } > \Delta\right] \\
\leq \exp\left\{- \frac{(\Delta/N) N e^{-B(T)T}}{CT}\right\}
\end{eqnarray*}
The announced convergence in law follows.

To solve the PDE, we set $\Psi(t,\theta)=e^{\theta_1^2 A(t)+
\theta_1\theta_2 C(t)+ \theta^2_2 B(t)}$. Then, substituting in the initial equation, and identifying the coefficients, we get the system $(\mathbf{M})$.$\  \  \otimes$

\paragraph{Proof of Theorem 2.4} We want to prove that the process $Z$ has the same law as the limit as $N\to \infty$ of the difference between the Langevin approximation linearized around the deterministic solutions and the deterministic solution itself, scaled by $\sqrt{N}$. We write it in the general case, not only in dimension two as above.
\noindent First we identify the equations satisfy by the moments of $Z$ starting from the equation satisfied by the characteristic function. We make the ersatz:
$$\psi(t,\theta)=e^{-\frac{1}{2}\theta \Gamma(t) \theta^T}$$
\noindent The matrix $\Gamma_t$ corresponds to the variance/covariance matrix.
We plug this expression into the equation satisfied by $\psi$ as
given in theorem 2.3:
\begin{eqnarray*}
\frac{\partial\Psi}{\partial t} &=& \sum_{j=1}^q \left\{\sum_{l\in L}\sum_{k=1}^{r_j}\theta^{(j)}_k\frac{\partial b_{j,k}}{\partial x_l}\frac{\partial\Psi}{\partial \theta_l} -\frac{1}{2}\sum_{k,l=1}^{r_j}\theta^{(j)}_{k}\theta^{(j)}_lG^{(j)}_{k,l}\Psi\right\}\\
&+&\sum_{m=1}^p\sum_{l\in L}\theta_m\frac{\partial f^m}{\partial
x_l}\frac{\partial\Psi}{\partial \theta_l}
\end{eqnarray*}
The ensemble of indices $L$ can be writen $L=L_v\cup L_u$ where
$L_v=\{1\leq m \leq p\}$ and $L_u=\{(j,k),\ 1\leq j\leq q,\ 1\leq
k\leq r_j\}$. To identify the equations satisfied by $\Gamma_{ab}$ we
distinguish the following cases:
\begin{itemize}
\item $a\in L_v$ and $b\in L_v$:
$\frac{1}{2}\Gamma'_{ab}=\displaystyle{\sum_{l\in L}} \left[ \frac{\partial f^a}{\partial x_l}\Gamma_{bl} + \frac{\partial f^b}{\partial x_l}\Gamma_{al} \right]$
\item $a\in L_v$ and $b\in L_u$, $b=(j,k)$:
$\frac{1}{2}\Gamma'_{ab}=\displaystyle{\sum_{l\in L}} \left[ \frac{\partial b_{j,k}}{\partial x_l}\Gamma_{al} + \frac{\partial f^a}{\partial x_l}\Gamma_{bl} \right]$
\item $a\in L_u$, $a=(j,k)$ and $b\in L_v$:
$\frac{1}{2}\Gamma'_{ab}=\displaystyle{\sum_{l\in L}} \left[ \frac{\partial b_{j,k}}{\partial x_l}\Gamma_{bl} + \frac{\partial f^b}{\partial x_l}\Gamma_{al} \right]$
\item $a\in L_u$, $a=(j,k)$ and $b\in L_v$, $b=(j',k')$:
$$\frac{1}{2}\Gamma'_{ab}=\displaystyle{\sum_{l\in L}} \left[ \frac{\partial b_{j,k}}{\partial x_l}\Gamma_{bl} + \frac{\partial b_{j',k'}}{\partial x_l}\Gamma_{al} \right]+\frac{1}{2}G_{k,k'}^{(j)}\mathbf{1}_{j=j'}$$
\end{itemize}

\noindent We then write the equations satisfied by $K^{(N)}(t)=\sqrt{N}(\tilde{X}_N(t)-X(t))=(Y_N^m, P_N^{j,k})$, where $\tilde{X}_N$ is the Langevin approximation defined in section 2.4, and where $X$ is the deterministic limit:

$$\begin{array}{ll}
dY_N^m=\sqrt{N}(f^m(\tilde{X}_N)-f^m(X))dt\\
dP_N^{j,k}=\sqrt{N}(b_{j,k}(\tilde{X}_N)-b_{j,k}(X))dt+\sigma^{(j)}(\tilde{X}_N)dW^j_t
\end{array}$$

When we linearize around the deterministic solution, we obtain the
following equations:

$$\begin{array}{ll}
d\tilde{Y}_N^m=\displaystyle{\sum_{l\in L}}\frac{\partial f^m}{\partial x_l}K^{(N)}_ldt\\
d\tilde{P}_N^{j,k}=\displaystyle{\sum_{l\in L}} \frac{\partial
b_{j,k}}{\partial
x_l}K^{(N)}_ldt+\left(\displaystyle{\sum_{k'=1}^{r_j}}\sigma^{(j)}_{k,k'}(X)+\frac{1}{\sqrt{N}}\Omega^{(j)}_{k,k'}\right)dW^{j,k'}_t
\end{array}$$

where the terms $\frac{1}{\sqrt{N}}\Omega^{(j)}_{k,k'}$ comes from
the linearization of $\sigma^{(j)}_{k,k'}(\tilde{X}_N)$, we do not need to
specify them here because they go to zero as $N\to\infty$.

It is now clear that the moments equations for this linear diffusion
system tends the system satisfied by $\Gamma_{ab}$ as
$N\to\infty$.

\paragraph{Proof of Theorem 2.5}
The convergence of $X_N$ to $X$ a.s. uniformly on finite time intervals, obtained in Theorem 2.1, implies that $\tau_N\to \tau$ a.s. In order apply Theorem 2.3, let us introduce $Z_N$ through the following decomposition:
\begin{eqnarray*}
\sqrt{N}(\phi(X(\tau))-\phi(X(\tau_N))) &=& \sqrt{N}\left[\phi \left(X(\tau_N)+\frac{1}{\sqrt{N}}Z_N(\tau_N)\right)-\phi \left(X(\tau_N)\right)\right]\\
&-& \sqrt{N}\phi(X_N(\tau_N))
\end{eqnarray*}
As $N\to \infty$, we claim that the right hand side converges in law to $\nabla \phi(X(\tau)).Z(\tau)$ since $\sqrt{N}\phi(X_N(\tau_N))$ converges in law to zero. Indeed,  as $\phi(X_N(\tau_N))\leq 0$ and $\phi(X_N(\tau_N^-))\geq 0$,
\begin{eqnarray*}
|\sqrt{N}\phi(X_N(\tau_N))|&\leq & |\sqrt{N}(\phi(X_N(\tau_N))-\phi(X_N(\tau_N^-))|
\end{eqnarray*}
There exists $\theta_N$ on the line between $X_N(\tau_N)$ and $X_N(\tau_N^-)$ such that $$|\sqrt{N}(\phi(X_N(\tau_N))-\phi(X_N(\tau_N^-))|=|\nabla \phi(\theta_N).(Z_N(\tau_N)-Z_N(\tau_N^-))|$$ which converges in law to zero since $Z_N\to Z$ and $Z$ is continuous. The claim follows.
By continuity, $\phi(X(\tau))=0$, so that $\sqrt{N}(\phi(X(\tau))-\phi(X(\tau_N)))$ is asymptotic to $$-\nabla \phi(X(\tau)).F(X(\tau))\sqrt{N}(\tau_N-\tau)$$ Thus $\sqrt{N}(\tau_N-\tau)$ converges in law to $\pi(\tau)$. To finish the proof we remark that $\sqrt{N}(X_N(\tau_N)-X(\tau))=Z_N(\tau_N)+\sqrt{N}(X(\tau_N)-X(\tau))$ which converges in law to $Z(\tau)+\pi(\tau)F(X(\tau))$.

\appendix

\section{Comparison between two deterministic limits of different stochastic Hodgkin-Huxley models}

We want to compare the two following systems deterministic (A.1) and (A.2), with
$f,\alpha,\beta$ continuously differentiable functions, $\alpha$ and $\beta$ non-negative, $k$ an integer $\geq 1$:
\begin{eqnarray}
\left\{
    \begin{array}{lll}
       \frac{dV}{dt}  =  f(V,u^k)\\
       \frac{du}{dt}  =  (1-u)\alpha(V)-u\beta(V)\\

    \end{array}
    \right.
\end{eqnarray}

\begin{eqnarray}
\left\{
    \begin{array}{llll}
       \frac{d\hat{V}}{dt}  =  f(\hat{V},x_k )\\

       \frac{dx_j}{dt}=(k-j+1) x_{j-1}\alpha(\hat{V})+(j+1)x_{j+1}\beta(\hat{V})\\
       -x_j\left(j\beta(\hat{V})+(k-j)\alpha(\hat{V})\right)\\
       \forall 0\leq j\leq k

    \end{array}
    \right.
\end{eqnarray}

System (A.1) corresponds to the classical ``Hodgkin-Huxley''
model, with only two variables for simplicity, and the system (A.2) is
a $(k+2)$-dimensional system, where $x_j,\ 0\leq j\leq k$ is the
proportion of channels in the state $j$, and $j=k$ is the
\textit{open} state.

\paragraph{Proposition}
Let $V_0 \in R$ and $u_0 \in [0,1]$. If the following conditions on the initial values are satisfied:\\
$V(0)=\hat{V}(0)=V_0$ and $\forall 0\leq j \leq k,\
C_k^ju(0)^{k-j}(1-u(0))^j=x_{k-j}(0)=C_k^ju_0^{k-j}(1-u_0)^j$
Then, for all $t\geq 0$, $V(t)=\hat{V}(t)$ (same potential) and $u(t)^k=x_k(t)$ (the proportion of open channels is $u(t)^k$).\\
Moreover, for all $1 \leq j \leq k$, for all $t\geq 0$, $x_{k-j}(t)=
C_k^ju(t)^{k-j}(1-u(t))^j$

\paragraph{Proof} Consider $(V,u)$ the unique solution of (1) for $V(0)=V_0$ and $u(0)=u_0$. Let $y_j(t)= C_k^ju(t)^{k-j}(1-u(t))^j$, $0 \leq j \leq k$. Then $(V,y_k,...y_0)$ is a solution of (2) (just need to compute $y_j'$ and write it in function of $y_{j-1}$ and $y_{j+1}$). As the initial values are equal (by hypothesis) : $x_{k-j}(0)=C_k^ju_0^{k-j}(1-u_0)^j=y_{k-j}(0)$, by uniqueness $(V,y_k,...,y_0)=(\hat{V},x_k,...x_0)$ for all $t\geq 0$.

\paragraph{Remark} The result is essentially the same for more complicated Markov schemes, as the sodium multistate Markov
model.

\section{Moments equations for linearized Langevin approximation}
From \textbf{Theorem 2.2}, one can build a diffusion approximation $(\tilde{V}_N,\tilde{e}_N)$ of the stochastic hybrid process $(V_N,e_N)$ given in the Example of section 2.1:\\

$$\left\{
\begin{array}{llll}
d\tilde{V}_N(t)=f(\tilde{V}_N(t),\tilde{e}_N(t))dt \\
d\tilde{e}_N(t)=b(\tilde{e}_N(t),\tilde{V}_N(t))dt+\sqrt{\frac{\lambda(\tilde{e}_N(t),\tilde{V}_N(t))}{N}}dB_t \\
b(v,e)=\left[(1-e)\alpha(v)-e\beta(v)\right]\\
\lambda(v,e)=\left[(1-e)\alpha(v)+e\beta(v)\right]
\end{array}
\right.$$
\noindent We want to write the moments equations for the linearized version of $$\left( \begin{array}{c} \tilde{P}_N \\ \tilde{Y}_N \end{array} \right) = \left( \begin{array}{c} \sqrt{N}\left(\tilde{e}_N-e\right) \\ \sqrt{N}\left(\tilde{V}_N-V\right) \end{array} \right)$$ with $(V,e)$ the deterministic solution. The linearized equations are given by:
$$\left\{
\right.
\begin{array}{ll}
dY^L_N=(f'_V Y^L_N +f'_e P^L_N)dt\\
dP^L_N=(b'_V Y^L_N +b'_e P^L_N)dt + \left[\sqrt{\lambda_t} +
\frac{1}{2\sqrt{N\lambda_t}}(\lambda'_V Y^L_N + \lambda'_e
P^L_N)\right]dB_t
\end{array}$$ with $\lambda_t=\lambda(V(t),e(t))$. We define $m^N_1=\mathbf{E}[Y^L_N]$, $m^N_2=\mathbf{E}[P^L_N]$, $S^N_1=\mathbf{E}[(Y^L_N-m_1)^2]$, $S^N_2=\mathbf{E}[(P^L_N-m_2)^2]$ and $C^N_{12}=\mathbf{E}[(Y^L_N-m_1)(P^L_N-m_2)]$.
Then we have the following system of 5 equations:
$$\left\{
\right.
\begin{array}{lll}
\frac{dm^N_1}{dt}&=&f'_V m_1 + f'_e m_2\\
\frac{dm^N_2}{dt}&=&b'_V m_1 + b'_e m_2\\

\end{array}$$
$$\left\{
\right.
\begin{array}{llll}
\frac{dS^N_1}{dt}&=&2f'_V S_1 + 2f'_e C_{12}\\
\frac{dS^N_2}{dt}&=&2b'_e S_2 + 2b'_V C_{12}+ \left[\sqrt{\lambda_t}+\frac{1}{2\sqrt{N \lambda}}(\lambda'_N m_1 + \lambda'_e m_2)\right]^2 \\
&+& \left(\frac{\lambda'_V}{2\sqrt{N\lambda_t}}\right)^2 S_1 + \left(\frac{\lambda'_e}{2\sqrt{N\lambda_t}}\right)^2 S_2 + 2 \frac{\lambda'_V \lambda'_e}{4N\lambda_t} C_{12}\\
\frac{dC^N_{12}}{dt}&=&b'_V S_1 + f'_e S_2 + (f'_V+b'_e)C_{12}

\end{array}$$
At the limit $N\to \infty$ and with $A=-2S_2$, $B=-2S_1$ and
$C=-C_{12}$ this system is the same as the one found in application
of \textbf{Theorem 2.3} in section 3.
\section{Moments equations for the Morris-Lecar system}
The moments equations used in section 3.5.1 and 3.5.2 are the following linear non-homogeneous system of differential equations:
$$\left[ \begin{array}{c} S_m \\ S_n \\ S_v \\ C_{mv} \\C_{nv} \\ C_{mn} \end{array} \right]=M(t)\left[ \begin{array}{c} S_m \\ S_n \\ S_n \\ C_{mv} \\C_{nv} \\ C_{mn} \end{array} \right]+\left[ \begin{array}{c} B_1 \\ B_2 \\ 0\\ 0 \\0 \\ 0 \end{array} \right]$$
with  $$M(t)=\left[ \begin{array}{cccccc}
2\frac{\partial F_m}{\partial m} & 0 & 0 & 2\frac{\partial F_m}{\partial V} & 0 & 0\\
0 & 2\frac{\partial F_n}{\partial n} & 0 & 0 & 2\frac{\partial F_n}{\partial V} & 0\\
0 & 0 & 2\frac{\partial F_v}{\partial V} & 2\frac{\partial F_v}{\partial m} & 02\frac{\partial F_v}{\partial n} & 0\\
\frac{\partial F_v}{\partial m} & 0 & \frac{\partial F_m}{\partial V} & \frac{\partial F_v}{\partial V}+\frac{\partial F_m}{\partial m} & 0 & \frac{\partial F_v}{\partial n}\\
0 & \frac{\partial F_v}{\partial n} & \frac{\partial F_n}{\partial V} & 0 & \frac{\partial F_v}{\partial V}+\frac{\partial F_n}{\partial n} & \frac{\partial F_v}{\partial m}\\
0 & 0 & 0 & \frac{\partial F_n}{\partial V} & \frac{\partial F_m}{\partial V} & \frac{\partial F_m}{\partial m}+\frac{\partial F_n}{\partial n}\\  
   \end{array} \right]$$ all the functions being evaluated at $X(t)=(V(t),m(t),n(t))$ solution of (3.5-3.7) and with $B_1(t)=(1-m(t))\alpha_m(V(t))+m(t)\beta_m(V(t))$, $B_2(t)=(1-n(t))\alpha_n(V(t))+n(t)\beta_n(V(t))$.

\section*{Aknowledgements}
During this work, G.Wainrib, supported by a fellowship from Ecole
Polytechnique, has been hosted by the Institut Jacques Monod and the
Laboratoire de Probabilit\'es et Mod\`eles Al\'eatoires, and wants
to thank both of them for their hospitality. This work has been supported by the project MANDy, ANR-09-BLAN-0008-01 of the Agence Nationale de la Recherche (ANR).

\bibliographystyle{plain}
\bibliography{biblio.bib}

\begin{thebibliography}{10}

\bibitem{Austin2007}
T.D. Austin.
\newblock {The emergence of the deterministic Hodgkin--Huxley equations as a
  limit from the underlying stochastic ion-channel mechanism}.
\newblock {\em Ann. Appl. Probab}, 18(4):1279--1325, 2008.

\bibitem{blom2006}
H.A.P. Blom and J.~Lygeros.
\newblock {Stochastic hybrid systems(theory and safety critical applications)}.
\newblock {\em Lecture notes in control and information sciences}, 2006.

\bibitem{chicone99}
C.C. Chicone.
\newblock {\em {Ordinary Differential Equations with Applications}}.
\newblock Springer, 1999.

\bibitem{cronin}
J.~Cronin.
\newblock {\em {Mathematical Aspects of Hodgkin-Huxley Neural Theory}}.
\newblock Cambridge University Press, 1987.

\bibitem{darling2005}
RWR Darling and JR~Norris.
\newblock {Structure of large random hypergraphs}.
\newblock {\em Ann. Appl. Probab}, 15(1A):125--152, 2005.

\bibitem{Davis1984}
M.~Davis.
\newblock Piecewise-deterministic markov processes: a general class of
  non-diffusion stochastic models.
\newblock {\em Journal of the royal statistical society (B)}, 43,3:353--388,
  1984.

\bibitem{defelice1993}
L.J. DeFelice and A.~Isaac.
\newblock {Chaotic states in a random world: Relationship between the nonlinear
  differential equations of excitability and the stochastic properties of ion
  channels}.
\newblock {\em Journal of Statistical Physics}, 70(1):339--354, 1993.

\bibitem{destexhe1994}
A.~Destexhe, Z.F. Mainen, and T.J. Sejnowski.
\newblock {Synthesis of models for excitable membranes, synaptic transmission
  and neuromodulation using a common kinetic formalism}.
\newblock {\em Journal of Computational Neuroscience}, 1(3):195--230, 1994.

\bibitem{pankratova}
A.~V.~Polovinkin E.~V.~Pankratova and E.~Mosekilde.
\newblock Resonant activation in a stochastic hodgkin-huxley model: Interplay
  between noise and suprathreshold driving effects.
\newblock {\em The European Physical Journal B - Condensed Matter and Complex
  Systems}, 45(3):391--397, 2005.

\bibitem{ethierkurtz}
S.N. Ethier and T.G. Kurtz.
\newblock {\em Markov Processes, Characterization and Convergence}.
\newblock John Wiley and Sons, Inc., 1986.

\bibitem{Faisal08}
AA. Faisal, LP. Selen, and DM. Wolpert.
\newblock Noise in the nervous system.
\newblock {\em Nat Rev Neurosci}, 9(4):292--303, April 2008.

\bibitem{Fox1994}
R.F. Fox and Y.~Lu.
\newblock {Emergent collective behavior in large numbers of globally coupled
  independently stochastic ion channels}.
\newblock {\em Physical Review E}, 49(4):3421--3431, 1994.

\bibitem{gollish}
T.~Gollisch and M.~Meister.
\newblock Rapid neural coding in the retina with relative spike latencies.
\newblock {\em Science}, 319(1):1108--1111, 2008.

\bibitem{guckoliva}
J.~Guckenheimer and R.A. Oliva.
\newblock Chaos in the hodgkin-huxley model.
\newblock {\em SIAM J. Appl. Dynam. Syst.}, 1:105--114, 2002.

\bibitem{hespanha2005}
J.P. Hespanha.
\newblock {A model for stochastic hybrid systems with application to
  communication networks}.
\newblock {\em Nonlinear Analysis}, 62(8):1353--1383, 2005.

\bibitem{HH52}
A.~L. Hodgkin and A.~F. Huxley.
\newblock A quantitative description of membrane current and its application to
  conduction and excitation in nerve.
\newblock {\em J Physiol}, 117(4):500--544, August 1952.

\bibitem{siam}
Kre\v{s}imir Josi\'{c} and Robert Rosenbaum.
\newblock Unstable solutions of nonautonomous linear differential equations.
\newblock {\em SIAM Rev.}, 50(3):570--584, 2008.

\bibitem{Kallenberg}
O.~Kallenberg.
\newblock {\em {Foundations of Modern Probability}}.
\newblock Springer, 1997.

\bibitem{keener}
J.~Keener.
\newblock Invariant manifold reductions for markovian ion channel dynamics.
\newblock {\em Journal of Mathematical Biology}, 58(3):447--57, 2009.

\bibitem{Kouretas2006}
P.~Kouretas, K.~Koutroumpas, J.~Lygeros, and Z.~Lygerou.
\newblock {Stochastic Hybrid Modeling of Biochemical Processes}.
\newblock {\em Stochastic Hybrid Systems}, 2006.

\bibitem{friesen}
H.~Krausz and W.O. Friesen.
\newblock The analysis of nonlinear synaptic transmission.
\newblock {\em Journal of General Physiology}, 70:243--265, 1977.

\bibitem{Kurtz1971}
T.G Kurtz.
\newblock Limit theorems for sequences of jump markov processes approximating
  ordinary differential processes.
\newblock {\em J. Appl Prob}, 8:344--356, 1971.

\bibitem{LygPNAS}
J.~Lygeros, K.~Koutroumpas, S.~Dimopoulos, I.~Legouras, P.~Kouretas,
  C.~Heichinger, P.~Nurse, and Z.~Lygerou.
\newblock {Stochastic hybrid modeling of DNA replication across a complete
  genome}.
\newblock {\em Proceedings of the National Academy of Science},
  105(34):12295--12300, 2008.

\bibitem{mainen1995msi}
ZF~Mainen, J.~Joerges, JR~Huguenard, and TJ~Sejnowski.
\newblock {A model of spike initiation in neocortical pyramidal neurons.}
\newblock {\em Neuron}, 15(6):1427--39, 1995.

\bibitem{morris1981}
C.~Morris and H.~Lecar.
\newblock {Voltage oscillations in the barnacle giant muscle fiber}.
\newblock {\em Biophysical Journal}, 35(1):193--213, 1981.

\bibitem{note}
K.~Pakdaman, M.~Thieullen, and G.~Wainrib.
\newblock A note on the large deviations and exit point for diffusion
  approximation of jump process : Markov vs. langevin.
\newblock {\em in preparation}.

\bibitem{refHH}
J.~Rinzel and R.~Miller.
\newblock Numerical calculation of stable and unstable periodic solutions to
  the hodgkin-huxley equations.
\newblock {\em Mathematical Bioscience}, 49:27--59, 1980.

\bibitem{rowat2007iis}
P.~Rowat.
\newblock {Interspike Interval Statistics in the Stochastic Hodgkin-Huxley
  Model: Coexistence of Gamma Frequency Bursts and Highly Irregular Firing}.
\newblock {\em Neural Computation}, 19(5):1215, 2007.

\bibitem{segundo94}
JP~Segundo, J.F. Vibert, K.~Pakdaman, M.~Stiber, and O.D. Martinez.
\newblock {Chapter 13 Noise and the Neurosciences: A Long History, a Recent
  Revival and Some Theory}.
\newblock {\em in Origins: Brain and Self Organization}, 1994.

\bibitem{shuai2003}
JW~Shuai and P.~Jung.
\newblock {Optimal ion channel clustering for intracellular calcium signaling}.
\newblock {\em Proceedings of the National Academy of Sciences},
  100(2):506--512, 2003.

\bibitem{skaugen1979fbs}
E.~Skaugen and L.~Walloe.
\newblock {Firing behaviour in a stochastic nerve membrane model based upon the
  Hodgkin-Huxley equations.}
\newblock {\em Acta Physiol Scand}, 107(4):343--63, 1979.

\bibitem{Steinmetz2000}
P.N. Steinmetz, A.~Manwani, C.~Koch, M.~London, and I.~Segev.
\newblock Subthreshold voltage noise due to channel fluctuations in active
  neuronal membranes.
\newblock {\em Journal of Computational Neuroscience}, 9(16):133--148, 2000.

\bibitem{touboul}
Jonathan Touboul and Olivier Faugeras.
\newblock First hitting time of double integral processes to curved boundaries.
\newblock {\em Advances in Applied Probability}, 40(2):501--528, 2008.

\bibitem{Tuckwell87}
H.~Tuckwell.
\newblock Diffusion approximations to channel noise.
\newblock {\em J. theor. Biol}, 127:427--438, 1987.

\bibitem{Vandenberg1991}
CA~Vandenberg and F.~Bezanilla.
\newblock {A sodium channel gating model based on single channel, macroscopic
  ionic, and gating currents in the squid giant axon}.
\newblock {\em Biophysical Journal}, 60(6):1511--1533, 1991.

\bibitem{thorpe}
R.~VanRullen, R.~Guyonneau, and S.~Thorpe.
\newblock Spike times make sense.
\newblock {\em TRENDS in Neurosciences}, 28(1):1--4, 2005.

\bibitem{white2000cnn}
J.A. White, J.T. Rubinstein, and A.R. Kay.
\newblock {Channel noise in neurons}.
\newblock {\em Trends in Neurosciences}, 23(3):131--137, 2000.

\end{thebibliography}

\end{document}